\documentclass[a4paper,reqno,11pt]{amsart}
\usepackage[T1]{fontenc}
\usepackage[utf8]{inputenc}
\usepackage{amssymb, amsmath, amsthm, graphicx, enumitem, color}
\usepackage[usenames,dvipsnames,svgnames,table]{xcolor}
\usepackage[foot]{amsaddr}
\usepackage{comment}

\usepackage[longnamesfirst,numbers,sort&compress]{natbib}

\makeatletter
\def\thm@space@setup{
\thm@preskip=4mm
\thm@postskip=0mm
}
\makeatother

\usepackage{hyperref}
\hypersetup{
    pdftitle={Quickly excluding an apex forest},
    pdfauthor={Jędrzej Hodor, Hoang La, Piotr Micek, Clément Rambaud},
    colorlinks,
    linkcolor={RoyalBlue},
    citecolor={RubineRed},
    urlcolor={blue!80!black}
}

\usepackage{cleveref}
\usepackage{mathtools}
\usepackage{thmtools, thm-restate}

\setenumerate{label=\textup{(\roman*)}, noitemsep, topsep=-\parskip, labelindent=.2em, leftmargin=*, widest=iii,}
\setitemize{noitemsep, topsep=-\parskip, labelindent=.2em, leftmargin=*, widest=iii,}

\DeclarePairedDelimiter\set{\{}{\}}

\theoremstyle{plain}
\newtheorem{thm}{Theorem}
\newtheorem*{thm*}{Theorem}
\newtheorem{theorem}[thm]{Theorem}

\newtheorem{lemma}[thm]{Lemma}
\newtheorem*{lemma*}{Lemma}
\newtheorem{cor}[thm]{Corollary}
\newtheorem*{cor*}{Corollary}

\newtheorem*{corollary*}{Corollary}
\newtheorem{obs}[thm]{Observation}

\theoremstyle{remark}
\newtheorem{problem}{Problem}

\newtheorem{claim}{Claim}[thm]
\newtheorem*{claim*}{Claim}
\crefname{obs}{Observation}{Observations}
\theoremstyle{definition}

\newtheorem*{conj*}{Conjecture}
\crefname{lem}{Lemma}{Lemmas}
\crefname{thm}{Theorem}{Theorems}
\crefname{cor}{Corollary}{Corollaries}

\newenvironment{proofclaim}[1][]
	{\begin{proof}[Proof of the claim] }{\end{proof}}

\newcommand{\td}{\operatorname{td}}
\newcommand{\pw}{\operatorname{pw}}
\newcommand{\tw}{\operatorname{tw}}
\newcommand{\ltd}{\operatorname{ltd}}
\newcommand{\lpw}{\operatorname{lpw}}

\newcommand{\gm}{f_{\boxplus}}

\newcommand{\Oh}{\mathcal{O}}

\newcommand{\cgW}{\mathcal{W}}

\newcommand{\dist}{\mathrm{dist}}

\let\leq\leqslant
\let\geq\geqslant

\let\subset\subseteq

\let\epsilon\varepsilon

\DeclareMathOperator\diam{diam}

\DeclareMathOperator\tn{tn}

\renewcommand{\setminus}{-}
\renewcommand{\vec}{\overrightarrow}

\title{Quickly excluding an apex-forest}

\begin{document}

\author[Hodor]{Jędrzej Hodor}
\address[J.~Hodor]{Theoretical Computer Science Department, 
Faculty of Mathematics and Computer Science and Doctoral School of Exact and Natural Sciences, Jagiellonian University, Krak\'ow, Poland}
\email{jedrzej.hodor@gmail.com}

\author[La]{Hoang La}
\address[H.~La]{LISN, Universit\'e Paris-Saclay, CNRS, Gif-sur-Yvette, France}
\email{hoang.la.research@gmail.com}

\author[Micek]{Piotr Micek}
\address[P.~Micek]{Theoretical Computer Science Department, 
Faculty of Mathematics and Computer Science, Jagiellonian University, Krak\'ow, Poland}
\email{piotr.micek@uj.edu.pl}

\author[Rambaud]{Clément Rambaud}
\address[C.~Rambaud]{Universit\'e C\^ote d'Azur, CNRS, Inria, I3S, Sophia-Antipolis, France}
\email{clement.rambaud@inria.fr}

\thanks{J.\ Hodor was supported by a Polish Ministry of Education and Science grant (Perły Nauki; PN/01/0265/2022).
P.\ Micek was supported by the National Science Center of Poland under grant UMO-2023/05/Y/ST6/00079 within the Weave-UNISONO program. 
C. Rambaud was supported by the French Agence Nationale de la Recherche under contract Digraphs ANR-19-CE48-0013-01.
An extended abstract of this paper has been accepted for publication in the proceedings of Eurocomb 2025.}

\begin{abstract}
We give a short proof that 
for every apex-forest $X$ on at least two vertices, graphs excluding $X$ as a minor have layered pathwidth at most $2|V(X)|-3$. This improves upon a result by Dujmovi\'c, Eppstein, Joret, Morin, and Wood (SIDMA, 2020). 
Our main tool is a structural result about graphs excluding a forest as a rooted minor, which is of independent interest.
We develop similar tools for treedepth and treewidth. 
We discuss implications for Erd\H{o}s-P\'osa properties of rooted models of minors in graphs.
\end{abstract}

\maketitle

\section{Introduction}

Within the seminal \emph{Graph minors} series, spanning from 1983 to 2010, 
Robertson and Seymour described the structure of graphs excluding a graph as a minor. 
One of many key insights of this series is the interplay between forbidding graphs as minors and treewidth or pathwidth.
Indeed, excluding a planar graph as a minor is equivalent to having bounded treewidth, which follows from the Grid Minor Theorem~\cite{GM5}. 
Similarly, excluding a forest as a minor is equivalent to having bounded pathwidth, which was proved in the first paper of the series~\cite{GM1}.
Another relevant statement following this pattern is that excluding a path as a minor is equivalent to having bounded treedepth; see e.g.~\cite[Chapter~6]{Sparsity-book}.

In this paper, we study analogous statements for excluding apex-type graphs as minors.
Recall that a graph is an \emph{apex graph} if it can be made planar by the removal of at most one vertex, and a graph is an \emph{apex-forest} if it can be made acyclic by the removal of at most one vertex.
It turns out that forbidding apex-type graphs as minors interplays with the layered versions of treewidth, pathwidth, and treedepth.
Dujmović, Morin, and Wood~\cite{layered-treewidth} proved that a minor-closed class of graphs excludes an apex graph if and only if it has bounded layered treewidth.
Similarly, Dujmović, Eppstein, Joret, Morin, and Wood in~\cite{layered-pathwidth}, proved that a minor-closed class of graphs excludes an apex-forest if and only if it has bounded layered pathwidth.
Our first contribution is a short and simple proof of the latter statement 
with an explicit and much better bound on layered pathwidth.
In what follows, for a graph $G$, we denote by $\tw(G)$, $\pw(G)$, $\td(G)$, and $\lpw(G)$ the treewidth, pathwidth, treedepth, and layered pathwidth of $G$ respectively.

\begin{theorem}\label{thm:lpw}
    For every apex-forest $X$ with at least two vertices, and for every graph $G$, if $G$ is $X$-minor-free, then 
    $\lpw(G)\leq 2|V(X)|-3$.
\end{theorem}

The main novelty in the proof of~\Cref{thm:lpw} is a version of path decomposition (and pathwidth) of a graph $G$ focused on some fixed subset $S$ of $V(G)$.
We denote the new parameter by $\pw(G,S)$, and we prove that for every forest $F$, if $G$ has no $S$-rooted model of $F$, then $\pw(G,S) \leq 2|V(F)|-2$, see~\Cref{thm:Spw}. 
This part of the proof follows closely the argument by Diestel~\cite{Diestel1995} showing that for every forest $F$, 
if $G$ has no $F$-minor, then $\pw(G)\leq |V(F)|-2$.

A graph is a \emph{fan} or (an \emph{apex-path}) if it becomes a path by the removal of at most one vertex. 
We introduce the concept of layered treedepth mimicking other layered parameters. It will be immediate that fans may have arbitrarily large layered treedepth.
Conversely, we prove that excluding a fan as a minor implies having bounded layered treedepth.
For a graph $G$, let $\ltd(G)$ be the layered treedepth of $G$. 

\begin{theorem}\label{thm:ltd}
    For every fan $X$ with at least three vertices, and for every graph $G$, if $G$ is $X$-minor-free, then 
    $\ltd(G)\leq \binom{|V(X)|-1}{2}$.
\end{theorem}

Similarly to the proof of \Cref{thm:lpw}, the proof of \Cref{thm:ltd} relies on a version of treedepth focused on some fixed subset $S$ of $V(G)$ that we denote by $\td(G,S)$.
The crucial property of this parameter is that for every path $P$, if $G$ has no $S$-rooted model of $P$, then $\td(G,S) \leq \binom{|V(P)|}{2}$, see~\Cref{thm:Std}. 

Following the definition of $\pw(G,S)$, one can also define a notion of treewidth focused on $S$, denoted by $\tw(G,S)$.
We show an approximate duality between this parameter and a version of tangles focused on $S$ proposed by Marx, Seymour, and Wollan~\cite{Marx2017}, see \Cref{thm:tw_and_tangle_numbers}. Combined with \Cref{thm:Marx_tangles}, a similar result to the main one of~\cite{Marx2017}, this yields a grid-minor theorem for this notion of treewidth, see \Cref{thm:Stw}.

The next two statements will follow immediately from the definitions of layered treedepth and pathwidth, \Cref{thm:ltd} and \Cref{thm:lpw}, respectively.
Recall that the \emph{diameter} of a graph $G$, denoted by $\diam(G)$, is the maximal distance between two vertices in $G$ taken over all pairs of vertices in $G$.

\begin{cor}\label{cor:td-diam}
    For every fan $X$ with at least two vertices, and for every connected graph $G$, if $G$ is $X$-minor-free, then $\td(G) \leq \binom{|V(X)|-1}{2}(\diam(G)+1)$.
\end{cor}

\begin{cor}\label{cor:pw-diam}
    For every apex-forest $X$ with at least two vertices, and for every connected graph $G$, if $G$ is $X$-minor-free, then $\pw(G) \leq (2|V(X)|-3)(\diam(G)+1)-1$. 
\end{cor}

\Cref{cor:td-diam,cor:pw-diam} are both optimal in the following sense.
There are fans of diameter $2$ and unbounded treedepth and there are apex-forests of diameter $2$ and unbounded pathwidth.
We also give a construction showing that the upper bound in~\Cref{cor:pw-diam} is tight up to a multiplicative constant, see \Cref{thm:td_rad_lowerbound}.

A natural strengthening of \Cref{thm:lpw} is the following (false) product structure statement: 
for every apex-forest $X$, there is a constant $c_X$ such that for every $X$-minor-free
graph $G$, we have $G \subseteq H \boxtimes P$\footnote{The \emph{strong product} $G_1\boxtimes G_2$ of two graphs $G_1$ and $G_2$ is the graph with vertex set $V(G_1\boxtimes G_2):=V(G_1)\times V(G_2)$ and that contains the edge with endpoints $(v,x)$ and $(w,y)$ if and only if
$vw\in E(G_1)$ and $x=y$; or $v=w$ and $xy\in E(G_2)$; or
$vw\in E(G_1)$ and $xy\in E(G_2)$.} for some graph $H$ with $\pw(H) \leq c_X$ and some path $P$.
The statement is false even for $X=K_3$, as Bose, Dujmovi\'c, Javarsineh, Morin, and Wood~\cite{Bose2022} proved that trees do not admit such a product structure.
Since $K_3$ is a fan, the construction in~\cite{Bose2022} also shows that the analogous strengthening of \Cref{thm:ltd} does not hold.

In \Cref{sec:prelimiaries}, we give all necessary definitions.
In \Cref{sec:rooted-minor}, we state our main technical contribution -- a collection of graph decompositions focused on a prescribed subset of vertices. 
In \Cref{sec:lpw}, we prove \Cref{thm:lpw} and the main properties of $\pw(G,S)$. 
In \Cref{sec:ltd}, we prove \Cref{thm:ltd} and the main properties of $\td(G,S)$. 
In \Cref{sec:tw}, we prove the main properties of $\tw(G,S)$.
Additionally, in \Cref{sec:lower-bounds}, we give a lower bound construction for \Cref{cor:td-diam} and \Cref{cor:pw-diam}.
In \Cref{sec:EP-property}, we discuss further applications of our results concerning Erd\H{o}s-P\'osa properties of rooted models of graphs.
In \Cref{sec:open}, we state a few open problems.
In \Cref{appendix:models,appendix}, we provide some complementary material included for completeness.

\section{Preliminaries}\label{sec:prelimiaries}

We allow a graph to be the null graph.
Moreover, the null graph is a path, a tree, and a grid. 
For a non-negative integer $\ell$, we denote by $P_{\ell}$ a graph that is a path on $\ell$ vertices.
For non-negative integers $k$ and $\ell$, as $k\times \ell$ \emph{grid}, we refer to the graph with vertex set $\{(i,j)\mid i \in [k],j\in [\ell]\}$ and vertices $(i,j)$ and $(i',j')$ are adjacent if and only if $|i'-i|+|j'-j|=1$.
We use $\boxplus_\ell$ to denote the $\ell \times \ell$ grid.
We say that a graph $G$ is a \emph{grid} if there exist non-negative integers $k$ and $\ell$ such that $G$ is the $k\times\ell$ grid.

Let $G$ be a graph. 
The \emph{neighborhood of a vertex} $v$ is the set of neighbors of $v$, denoted by $N_G(v)$. 
The \emph{closed neighborhood} of $v$ is then $N_G[v] = N_G(v) \cup \{v\}$.
The \emph{neighborhood of a set} $R\subseteq V(G)$ is $N_G(R)=\bigcup_{v\in R}N_G(v)-R$. We drop the subscript when the graph is clear from context.

Let $H$ be a graph.
A \emph{model} of $H$ in $G$ is a family $(B_x \mid x \in V(H))$ of pairwise disjoint subsets of $V(G)$ such that: 
\begin{enumerate}
    \item for every $x \in V(H)$, the subgraph of $G$ induced by $B_x$ is non-empty and connected, and
    \item for every edge $xy \in E(H)$, there is an edge between $B_x$ and $B_y$ in $G$.
\end{enumerate}
The set $B_x$ for $x \in V(H)$ is called the \emph{branch set} of $x$ in the model. 
We say that $H$ is a \emph{minor} of $G$ if $G$ contains a model of $H$. 
Otherwise, we say that $G$ is \emph{$H$-minor-free}.

Let $S \subseteq V(G)$.
A model $(B_x \mid x\in V(H))$ of $H$ in $G$ is \emph{$S$-rooted} if $B_x\cap S\neq\emptyset$ for each $x \in V(H)$. 
Moreover, if $H$ is a plane graph, we say that 
a model $(B_x \mid x\in V(H))$ of $H$ in $G$ is \emph{$S$-outer-rooted} if $B_x\cap S\neq\emptyset$ for each vertex $x$ in the outer face of $H$.

A \emph{rooted forest} is a disjoint union of rooted trees.
The \emph{vertex-height} of a rooted forest $F$ is the maximum number of vertices on a path from a root to a leaf in $F$, and the \emph{depth} of a vertex $u \in V(F)$ is the number of vertices in the path between $u$ and the root of its component.
For two vertices $u$, $v$ in a rooted forest $F$, we say that $u$ is a \emph{descendant} of $v$ in $F$
if $v$ lies on the path from a root to $u$ in $F$.  
The \emph{closure} of $F$ is the graph with vertex set $V(F)$ and edge set $\set{vw\mid v\neq w\text{ and } v \text{ is a descendant of $w$ in $F$}}$. 
We say that $F$ is an \emph{elimination forest} of $G$ if $V(F)=V(G)$ and $G$ is a subgraph of the closure of $F$.
The \emph{treedepth} of a graph $G$, denoted by $\td(G)$, is $0$ if $G$ is empty, and otherwise is the minimum vertex-height of an elimination forest of $G$.

A \emph{tree decomposition} of $G$ is a pair $\mathcal{B} = \big(T,(W_x \mid x\in V(T))\big)$, 
where $T$ is a non-null tree and the sets $W_x$ for each $x \in V(T)$ are subsets of $V(G)$ 
called \emph{bags} satisfying:
\begin{enumerate}
\item for each edge $uv\in E(G)$ there is a bag containing both $u$ and $v$, and
\item for each vertex $v\in V(G)$ the set of vertices $x\in V(T)$ with 
$v\in W_x$  induces a non-empty subtree of $T$.
\end{enumerate}
The {\em width} of $\mathcal{B}$ is $\max\{|W_x|-1 \mid x\in V(T) \}$.
The \emph{treewidth} of $G$, denoted $\tw(G)$, is the minimum width of a tree decomposition of $G$.
A \emph{path decomposition} of $G$ is a tree decomposition $\big(T,(W_x \mid x\in V(T))\big)$, where $T$ is a path.
In that case, we simply write $(W_i \mid i \in [m])$ for $\big(T,(W_x \mid x\in V(T))\big)$, where $m=|V(T)|$, simply identifying $T$ with a path on $[m]$.
The \emph{pathwidth} of a graph $G$, denoted $\pw(G)$, is the minimum width of a path decomposition of $G$.

A \emph{layering} of $G$ is a sequence $(L_i \mid i\geq 0)$ of disjoint subsets of $V(G)$ whose union is $V(G)$ and such that for every $uv \in E(G)$ there is a non-negative integer $i$ such that $u,v \in L_i \cup L_{i+1}$.

Let $\mathcal{B}=\big(T,(B_x \mid x\in V(T))\big)$ be a tree decomposition of $G$ and let $\mathcal{L}=(L_i \mid i\geq 0)$. 
The \emph{width} of $(\mathcal{B},\mathcal{L})$ is $\max\{|B_x \cap L_i| \mid x\in V(T),\  i\geq0 \}$. 
The \emph{layered treewidth} of $G$ is the minimum width of a pair $(\mathcal{B},\mathcal{L})$, where $\mathcal{B}$ is a tree decomposition of $G$ and $\mathcal{L}$ is a layering of $G$. 
The \emph{layered pathwidth} of $G$ is 
the minimum width of a pair $(\mathcal{B},\mathcal{L})$, where $\mathcal{B}$ is a path decomposition of $G$ and $\mathcal{L}$ is a layering of $G$.

We propose a natural counterpart of the definitions above for treedepth. 
Let $F$ be an elimination forest of $G$, and let $\mathcal{L}=(L_i \mid i\geq 0)$ be a layering of $G$.
The \emph{width} of $(F,\mathcal{L})$ is $\max\{|R \cap L_i| \mid \text{$R$ is a root-to-leaf path in $F$}, i \geq 0\}$.
The \emph{layered treedepth} of $G$ is the minimum width of a pair $(F,\mathcal{L})$, where $F$ is an elimination forest of $G$ and $\mathcal{L}$ is a layering of~$G$. 

Let $G,H$ be two graphs. We denote by $G \cup H$ the graph with vertex set $V(G) \cup V(H)$ and edge set $E(G) \cup E(H)$.
Similarly, for every set $F$ of pairs of vertices of $G$, we denote by $G \cup F$ the graph with vertex set $V(G)$ and edge set $E(G) \cup F$.

A \emph{separation} of $G$ is a pair $(A, B)$ of subgraphs of $G$ such that $A\cup B = G$,
$E(A \cap B) = \emptyset$.
The \emph{order} of $(A,B)$ is $|V(A) \cap V(B)|$.
For $X,Y \subset V(G)$, an \emph{$X$--$Y$ path} is a path in $G$ that is either a one-vertex path with the vertex in $X \cap Y$ or a path with one endpoint in $X$ and the other endpoint in $Y$ such that no internal vertices are in $X \cup Y$.
We need the following well-known theorem.

\begin{theorem}[Menger's Theorem]
    Let $G$ be a graph and $X,Y \subset V(G)$.
    There exists a separation $(A,B)$ of $G$ such that $X \subset V(A)$, $Y \subset V(B)$, and there exists $|V(A) \cap V(B)|$ pairwise disjoint $X$--$Y$ paths.
\end{theorem}

\section{Excluding a rooted minor}\label{sec:rooted-minor}

In this section, we introduce a new family of graph parameters and we state some of their properties that are the key technical ingredients of the proofs of our main results.
However, we believe that the parameters and their properties are of independent interest.
We start with the new version of treedepth, and afterward, we discuss the new versions of pathwidth and treewidth.

Let $G$ be a graph and let $S\subseteq V(G)$.
An \emph{elimination forest} of $(G,S)$ is an elimination forest $F$ of $H$, an induced subgraph of $G$ such that $S$ is contained in $V(H)$ and for every component $C$ of $G-V(H)$,
either $V(C)$ has no neighbors in $G$, or
there is a root-to-leaf path in $F$ containing all the neighbors of $V(C)$ in $G$.
The \emph{treedepth} of $(G,S)$, denoted by $\td(G,S)$, is the minimum vertex-height of an elimination forest of $(G,S)$.
This notion is similar in its definition to the elimination distance to a given subgraph-closed class of graph $\mathcal{C}$ (see e.g.~\cite{paul2023universal}), which is defined as $\min \{\td(G,S) \mid \text{$S \subset V(G)$ such that every component of $G-S$ is in $\mathcal{C}$}\}$.

Recall that if a graph $G$ has no model of $P_{\ell}$, then $\td(G)<\ell$.
We prove an analogous result within the setting of $S$-rooted models of paths.

\begin{theorem}\label{thm:Std}
    For every positive integer $\ell$, for every graph $G$, and for every $S\subseteq V(G)$, if $G$ has no $S$-rooted model of $P_\ell$, then $\td(G,S)\leq \binom{\ell}{2}$. 
\end{theorem}

As already mentioned, \Cref{thm:Std} is the main ingredient of the proof of \Cref{thm:ltd}.
Actually, the intuition standing behind this is very simple.
For a vertex $u$ in a graph $G$, we set $S = N(u)$.
Now, if $G-u$ has a $S$-rooted model of a path $P_\ell$, then $G$ has a model of $P_\ell$ with a universal vertex added, and so $G$ has a model of every fan on $\ell+1$ vertices.

\begin{figure}[!htbp]
    \centering 
    \includegraphics[scale=1]{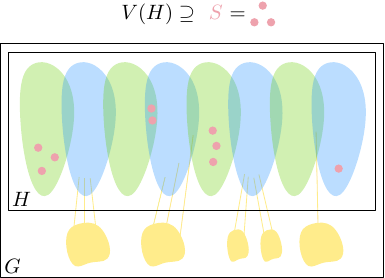} 
    \caption{The green and blue bags depict a path decomposition of $H$, an induced subgraph of $G$ such that $S \subset V(H)$.
    Each component of $G - V(H)$ (yellow) has all the neighbors in one of the bags. 
    Such a bag does not have to be unique.} \label{fig:pd_of_GS}
\end{figure}

Let $G$ be a graph and let $S\subseteq V(G)$.
A \emph{tree decomposition} (resp.\ \emph{path decomposition}) of $(G,S)$ is a tree decomposition (resp.\ path decomposition) $\mathcal{B}$ of $H$, an induced subgraph of $G$ such that $S$ is contained in $V(H)$,
and for every component $C$ of $G-V(H)$, 
there exists a bag of $\mathcal{B}$ containing all the neighbors of $V(C)$ in $G$. 
The \emph{treewidth} (resp.\ \emph{pathwidth}) of $(G,S)$, denoted by $\tw(G,S)$ (resp.\ $\pw(G,S)$), is the minimum width of a tree decomposition (resp.\ path decomposition) of $(G,S)$.
We illustrate the notion of a path decomposition of $(G,S)$ in \Cref{fig:pd_of_GS}.

In the first paper of the Graph Minors series~\cite{GM1}, Robertson and Seymour proved that if a graph $G$ has no model of a forest $F$, then $\pw(G)$ is bounded by a function of $|V(F)|$.
Bienstock, Robertson, Seymour, and Thomas~\cite{quickly-excluding-a-forest} obtained a tight result that if a graph $G$ has no model of a forest $F$, then $\pw(G)\leq |V(F)|-2$.
The most relevant work for our purposes is a beautiful and short proof of the inequality above given by Diestel in~\cite{Diestel1995}.
We prove an analogous result within the setting of $S$-rooted models of forests.

\begin{theorem}\label{thm:Spw}
    For every forest $F$ with at least one vertex, for every graph $G$, and for every $S\subseteq V(G)$, if $G$ has no $S$-rooted model of $F$, then $\pw(G,S)\leq 2|V(F)|-2$.
\end{theorem}

The Grid Minor Theorem can be generalized to the setting of $S$-outer-rooted models as follows.
In the following part of this section, let $\gm$ be the minimum function such that for every positive integer $\ell$, if a graph $G$ has no model of $\boxplus_\ell$, then $\tw(G) \leq \gm(\ell)$.
This function exists due to the Grid Minor Theorem and the best known upper bound is $\gm(\ell) \leq \ell^{9+o(1)}$ due to Chuzhoy and Tan~\cite{Chuzhoy2021}.

\begin{theorem}\label{thm:Stw}
    For every plane graph $H$, for every graph $G$, and for every $S\subseteq V(G)$, if $G$ has no $S$-outer-rooted model of $H$, then $\tw(G,S)\leq 3\gm(98304 \cdot |V(H)|^4) + 1$.
\end{theorem}

In this result, one cannot replace ``$S$-outer-rooted model'' with ``$S$-rooted model''.
Indeed, for every non-negative integer $\ell$, the graph $\boxplus_\ell$ with $S_\ell$ being the vertex set of the outer face has no $S_\ell$-rooted model of $K_4$, while $\tw(\boxplus_\ell,S_\ell) \geq \ell-1$.
The latter inequality follows from \Cref{lemma:EP_Stw} applied to the family of all the connected subgraphs that are the union of a row with a column.

It is well-known that every planar graph $H$ is a minor of $\boxplus_{2|V(H)|}$~\cite[statement 1.5]{Robertson1994}.
Following the same proof ideas, one can show that for every plane graph $H$, if a graph $G$ contains an $S$-outer-rooted model of $\boxplus_{2|V(H)|}$, then $G$ contains an $S$-outer-rooted model of $H$.
For completeness, we prove this in~\Cref{appendix:models}.
It follows that~\Cref{thm:Stw} is a consequence of the following more precise statement.
\begin{theorem}\label{thm:Stw-grids}
    For every positive integer $\ell$, for every graph $G$, and for every $S\subseteq V(G)$, if $G$ has no $S$-outer-rooted model of $\boxplus_\ell$, then $\tw(G,S)\leq 3\gm(6144\ell^4) + 1$.
\end{theorem}

We obtain \Cref{thm:Stw-grids} via tangles.
First, let us recall the definition of tangles in graphs.
Let $G$ be a graph and let $k$ be a positive integer.
Let $\mathcal{T}$ be a family of separations of $G$ of order less than $k$ in $G$. 
$\mathcal{T}$ is a \emph{tangle} of order $k$ in $G$ if
\begin{enumerate}[topsep=5pt-\parskip,label={\normalfont (T\arabic*)}]
    \item for every separation $(A,B)$ of order at most $k-1$ in $G$, $(A,B) \in \mathcal{T}$ or $(B,A) \in \mathcal{T}$, \label{item:T1}
    \item for every $(A_1,B_1),(A_2,B_2),(A_3,B_3) \in  \mathcal{T}$,
        $A_1 \cup A_2 \cup A_3 \neq G$, and \label{item:T2}
    \item for every $(A,B) \in \mathcal{T}$, $V(A) \neq V(G)$. \label{item:T3}
\end{enumerate}
Marx, Seymour, and Wollan~\cite{Marx2017} proposed a variant of tangles that is focused on a prescribed set.
Additionally, for a fixed $S\subseteq V(G)$, $\mathcal{T}$ is a \emph{tangle} of $(G,S)$ if it is a tangle of $G$ and 
\begin{enumerate}[resume*]
    \item \label{item:T4} for every $(A,B) \in \mathcal{T}$, $S \not\subseteq V(A)$.
\end{enumerate}
The \emph{tangle number} of $(G,S)$, denoted by $\tn(G,S)$, is the maximum order of a tangle of $(G,S)$. 
When $S = V(G)$, \cref{item:T4} is vacuous, and $\tn(G) = \tn(G,V(G))$ is the classical tangle number of $G$.
One of the cornerstones of structural graph theory is that the following graph parameters are functionally tied to each other: treewidth, tangle number, and the maximum integer $\ell$ such that a graph contains a model of~$\boxplus_\ell$.
An analog in the ``focused'' setting also holds.
\begin{theorem}[based on \cite{Marx2017}]\label{thm:Marx_tangles}
    For every positive integer $\ell$, for every graph $G$, and for every $S\subseteq V(G)$, if there is a tangle $\mathcal{T}$ of $(G,S)$ of order $3\gm(6144\ell^4) + 1$,
    then $G$ contains an $S$-outer-rooted model of $\boxplus_\ell$.
\end{theorem}

\Cref{thm:Marx_tangles} is a slight alteration of a result present in a paper of Marx, Seymour, and Wollan~\cite[statement~1.3]{Marx2017}.
What we need is not explicitly proved there, so for completeness we give a proof of~\Cref{thm:Marx_tangles} in~\Cref{appendix}.
To obtain~\Cref{thm:Stw}, we also functionally tie $\tn(G,S)$ and $\tw(G,S)$.

\begin{theorem}\label{thm:tw_and_tangle_numbers}
    For every graph $G$ with at least one vertex, and for every $S \subseteq V(G)$,
    \[
    \tn(G,S) - 1 \leq \tw(G,S) \leq 10\max\{\tn(G,S),2\}-12.
    \]
\end{theorem}
\Cref{thm:Marx_tangles} and \Cref{thm:tw_and_tangle_numbers} immediately imply \Cref{thm:Stw-grids}.

\section{Layered pathwidth}\label{sec:lpw}

The proof of \Cref{thm:Spw} follows Diestel's proof~\cite{Diestel1995} that graphs excluding a forest $F$ as a minor have pathwidth at most $|V(F)|-2$.
The notation follows a recent paper by Seymour~\cite{seymour2023shorter}.

Let $G$ be a graph, let $w$ be a positive integer, and let $S \subseteq V(G)$.
A separation $(A,B)$ is \emph{$(w,S)$-good} if it is of order at most $w$ and $(A,S\cap V(A))$ has a path decomposition of width at most $2w-2$ whose last bag contains $V(A) \cap V(B)$ as a subset.
When $(A,B)$, $(A',B')$ are separations of $G$, we write $(A,B) \leq (A',B')$, if $A \subseteq A'$ and $B \supseteq B'$.
If moreover
the order of $(A',B')$ is at most the order of $(A,B)$,
then we say that $(A',B')$ \emph{extends} $(A,B)$.
A separation $(A,B)$ in $G$ is \emph{maximal $(w,S)$-good} if it is $(w,S)$-good and for every
$(w,S)$-good separation $(A',B')$ extending $(A,B)$, we have $A'=A$ and $B'=B$.

\begin{figure}[!htbp]
    \centering 
    \includegraphics[scale=1]{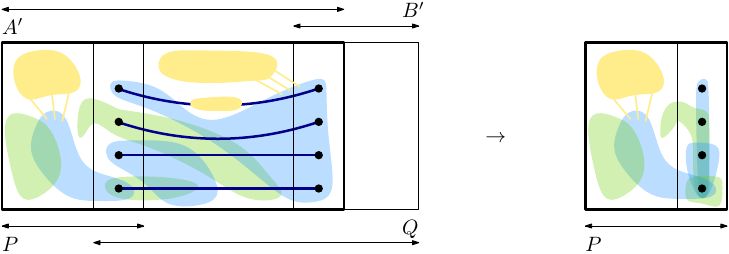} 
    \caption{
    Illustration of the proof of \Cref{lemma:key_lemma_wgood}.
    On the left, we depict the initial situation and on the right, we depict the result of applying the procedure from the lemma.
    Bags of path decompositions (of $(A',S \cap V(A'))$ on the left and of $(P,S \cap V(P))$ on the right) are green and blue alternately.
    In yellow, we show the components that are left after removing all vertices of respective path decompositions.
    The latter is obtained from the former by contracting $|V(P) \cap V(Q)|$ disjoint $V(P)$--$V(B')$ paths (in blue).
    }
    \label{fig:key_lemma_wgood}
\end{figure}

We start with a simple lemma illustrated in \Cref{fig:key_lemma_wgood}.

\begin{lemma}\label{lemma:key_lemma_wgood}
    Let $G$ be a graph, let $w$ be a positive integer, let $S \subseteq V(G)$, and let $(A',B')$ and $(P,Q)$ be two separations of $G$ with $(P,Q) \leq (A',B')$.
    If $(A',B')$ is $(w,S)$-good and there are $|V(P) \cap V(Q)|$ vertex-disjoint $V(P)$--$V(B')$ paths in $G$, then $(P,Q)$ is $(w,S)$-good.
\end{lemma}

\begin{proof}
    Suppose that $(A',B')$ is $(w,S)$-good and there are $|V(P) \cap V(Q)|$ vertex-disjoint $V(P)$--$V(B')$ paths $(R_x \mid x \in V(P) \cap V(Q))$ in $G$.
    Assume that $x \in V(R_x)$ for every $x \in V(P) \cap V(Q)$.
    Let $(W_i \mid i \in [m])$ be a path decomposition of $(A',S\cap V(A'))$ of width at most $2w-2$ with $V(A') \cap V(B') \subset W_m$.
    Let $(V_i \mid i \in [m])$ be obtained from $(W_i \mid i \in [m])$ by contracting  $R_x$ into a single vertex that we identify with $x$, for every $x \in V(P) \cap V(Q)$.
    In other words, $V_i = (W_i \cap V(P)) \cup \{x \in V(P) \cap V(Q) \mid V(R_x) \cap W_i \neq \emptyset\}$ for every $i \in [m]$.
    Observe that $|V_i| \leq |W_i|$ for every $i \in [m]$.
    We claim that $(V_i \mid i\in [m])$ is a path decomposition of $(P,S\cap V(P))$.

Let $X=V(P)\cap V(Q)$,
$H'=A'[\bigcup_{i\in[m]}W_i]$, and
$H=P[\bigcup_{i\in[m]}V_i]$.
Let $\mathcal D$ be the set of components of $A'-V(H')$.
For each $D\in\mathcal D$, choose $j_D\in[m]$ such that
$N_{A'}(V(D))\subseteq W_{j_D}$.
For each $i\in[m]$, define
\[
  \widehat W_i
  =W_i\cup
    \bigcup_{D\in\mathcal D, j_D=i}V(D).
\]
Then $(\widehat W_i\mid i\in[m])$ is a path decomposition of $A'$.
For every $x\in X$, the path $R_x$ meets $W_m$.
Thus, for every $D\in\mathcal D$ meeting $R_x$, the path $R_x$
also meets $N_{A'}(V(D))\subseteq W_{j_D}$.
Consequently, for all $x\in X$ and $i\in[m]$,
\[
  V(R_x)\cap\widehat W_i\ne\emptyset
  \quad\text{if and only if}\quad
  V(R_x)\cap W_i\ne\emptyset.
\]
Contracting each $R_x$ to $x$ and restricting the bags to $V(P)$
gives a path decomposition of $P$ with bags
\[
  \widehat V_i
  =(\widehat W_i\cap V(P))
   \cup\{x\in X\mid V(R_x)\cap\widehat W_i\ne\emptyset\},
  \text{ for } i\in[m].
\]
Indeed, the bags meeting each path form an interval, and
$V(R_x)\cap V(P)=\{x\}$ for every $x\in X$, so every edge
of $P$ remains covered.
Since $V(H)=(V(H')\cap V(P))\cup X$ and $\widehat V_i\cap V(H)=V_i$ for every $i\in[m]$,
the sequence $(V_i\mid i\in[m])$ is a path decomposition of $H$.
Also, $S\cap V(P)\subseteq V(H')\cap V(P)\subseteq V(H)$.

For each component $C$ of $P-V(H)$, let $D(C)\in\mathcal D$
be the component containing $C$.
Every vertex of $C$ belongs to exactly one bag $\widehat V_i$,
namely $\widehat V_{j_{D(C)}}$.
Since $N_P(V(C))\subseteq V(H)$, we have
\[
  N_P(V(C))
  \subseteq \widehat V_{j_{D(C)}}\cap V(H)
  =V_{j_{D(C)}}.
\]
    

    Finally, since $V(A') \cap V(B')\subseteq W_m$, we have $V(P) \cap V(Q) \subset V_m$.
    Additionally, $|V(P) \cap V(Q)| = |\{R_x \mid x \in V(P) \cap V(Q)\}| \leq |V(A') \cap V(B')| \leq w$.
    All this proves that $(P,Q)$ is $(w,S)$-good.
\end{proof}

We will use the following version of Menger's theorem in the proof of \Cref{lemma:find_wTS_spanning_separation}.
\begin{lemma}\label{lemma:Menger}
    Let $G$ be a graph and let $(A,B)$ and $(A',B')$ be two separations of $G$.
    If $(A,B) \leq (A',B')$, then there is a separation $(P,Q)$ of $G$ such that
    \begin{enumerate}
        \item $(A,B) \leq (P,Q) \leq (A',B')$, and 
        \item there are $|V(P) \cap V(Q)|$ pairwise disjoint $V(A)$--$V(B')$ paths in $G$.
    \end{enumerate}
\end{lemma}

\begin{lemma}\label{lemma:find_wTS_spanning_separation}
    Let $w$ be a positive integer, let $G$ be a graph, and let $S \subseteq V(G)$ such that $\pw(G,S) > 2w-2$.
    If $F$ is a forest on at most $w$ vertices, then there is a separation $(A,B)$ of $G$ such that
    \begin{enumerate}[label={\normalfont (s\arabic*)}]
        \item $|V(A) \cap V(B)| \leq |V(F)|$, \label{item:find_wTS_i}
        \item there is a $(V(A) \cap V(B))$-rooted model of $F$ in $A$, and \label{item:find_wTS_ii}
        \item $(A,B)$ is maximal $(w,S)$-good. \label{item:find_wTS_iii}
    \end{enumerate}
\end{lemma}

\Cref{item:find_wTS_ii} in the statement above implies $|V(A)\cap V(B)|\geq |V(F)|$, so by \cref{item:find_wTS_i} we have $|V(A)\cap V(B)| = |V(F)|$.

\begin{figure}[!htbp]
    \centering 
    \includegraphics[scale=1]{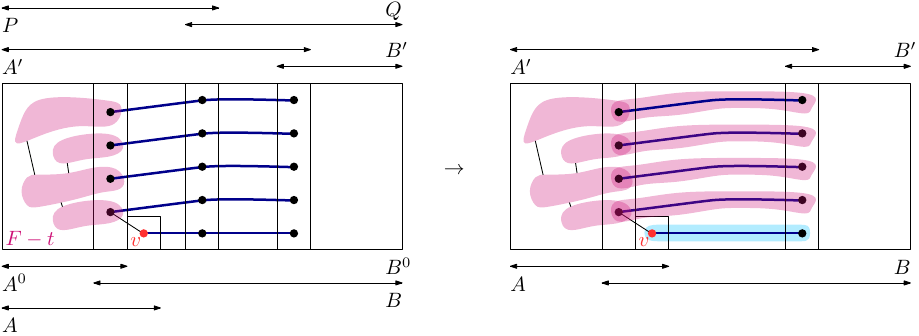} 
    \caption{An illustration of the proof of \Cref{lemma:find_wTS_spanning_separation}.
    We consider $F$ to be a forest ($|V(F)| = 5$ in the figure).
    In pink, we depict the branch sets of the rooted model of $F - t$.
    We argue that if $|V(P) \cap V(Q)| < |V(F)|$, then $(P,Q)$ contradicts the maximality of $(A^0,B^0)$.
    Hence, $V(A)$ is connected with $V(B')$ by $5$ pairwise disjoint paths.
    We add the blue branch set containing $v$ to the model and extend pink branch sets using the paths obtaining a $(V(A') \cap V(B'))$-rooted model of $F$ in $A'$.
    }
    \label{fig:maximal-sep}
\end{figure}

\begin{proof}
    We proceed by induction on $|V(F)|$.
    Suppose that $F$ is the null graph. 
    Since $(\emptyset,G)$ is $(w,S)$-good, then a maximal $(w,S)$-good separation $(A,B)$ extending $(\emptyset, G)$ satisfies~\ref{item:find_wTS_i}-\ref{item:find_wTS_iii}.
    Next, let $F$ be a non-empty forest on at most $w$ vertices.
    Let $t$ be a vertex of $F$ of degree at most one.

    By induction hypothesis for $F-t$, $G$ has a separation $(A^0,B^0)$ satisfying~\ref{item:find_wTS_i}-\ref{item:find_wTS_iii} for $F-t$.
    Let $(W_i \mid i \in [m])$ be a path decomposition of $(A^0,S\cap V(A^0))$ of width at most $2w-2$ with $V(A^0) \cap V(B^0) \subseteq W_m$.
    If $V(A^0)=V(G)$, then $(W_i \mid i\in [m])$ is a path decomposition of $(G,S)$,
    which contradicts the fact that $\pw(G,S) > 2w-2$.
    Hence $V(B^0) \setminus V(A^0) \neq \emptyset$.
    Let $(B_x \mid  x \in V(F-t))$ be a $(V(A^0) \cap V(B^0))$-rooted model of $F-t$ in $A^0$.
    If $t$ has degree $0$ in $F$, then choose a vertex $v \in V(B^0) \setminus V(A^0)$ arbitrarily.
    Otherwise, $t$ has a unique neighbor $s$ in $F$.
    By \ref{item:find_wTS_ii}, there is a vertex $u$ in $B_s \cap V(A^0) \cap V(B^0)$, and choose $v$ to be a neighbor of $u$ in $V(B^0) \setminus V(A^0)$.
    Such a neighbor exists as otherwise $(A^0 \cup \{uu' \mid uu' \in E(B^0)\},B^0 - u)$ is $(w,S)$-good, which contradicts the maximality of $(A^0,B^0)$.

    Let $(A,B)$ be the separation of $G$ defined by $A=G[V(A^0) \cup \{v\}]$ and $B=G[V(B^0)] \setminus E(A)$.
    Since $V(A^0) \cap V(B^0) \subset W_m$, and the neighborhood of $v$ in $A$ is contained in $V(A^0) \cap V(B^0)$, $(W_1, \dots, W_{m-1},W_m, V(A) \cap V(B))$ is a path decomposition of $(A,S\cap V(A))$.
    Moreover, since $|V(A)\cap V(B)| \leq |V(F)|\leq w \leq 2w-1$, this path decomposition is of order at most $2w-2$. 
    Therefore, $(A,B)$ is $(w,S)$-good, and so, there is a maximal $(w,S)$-good separation $(A',B')$ extending $(A,B)$ in $G$.
    In particular, $|V(A')\cap V(B')|\leq |V(F)|$.

    The next step of the proof is illustrated in \Cref{fig:maximal-sep}.
    By \Cref{lemma:Menger}, there exists a separation $(P,Q)$ such that $(A,B) \leq (P,Q) \leq (A',B')$ and there is a family $\mathcal{L}$ of $|V(P) \cap V(Q)|$ disjoint $V(A)$--$V(B')$ paths in $G$.
    If $|V(P) \cap V(Q)| \leq |V(F)|-1$, then by \Cref{lemma:key_lemma_wgood}, 
    since $(A',B')$ is $(w,S)$-good, 
    $(P,Q)$ is $(w,S)$-good as well.
    Since $(P,Q)$ extends $(A^0,B^0)$, and $v \in V(P) \setminus V(A^0)$, this contradicts the maximality of $(A^0,B^0)$.
    Hence $|V(P) \cap V(Q)| \geq |V(F)|$.
    Setting $B_t=\{v\}$ gives a $(V(A) \cap V(B))$-rooted model $(B_x \mid x \in V(F))$ of $F$ in $A$.
    Since $(A,B) \leq (A',B')$, every $V(A)$--$V(B')$ path is a $(V(A) \cap V(B))$--$(V(A')\cap V(B'))$ path contained in $V(B) \cap V(A')$.
    Therefore, the model can be extended using $|V(F)|$ paths in $\mathcal{L}$ yielding a $(V(A') \cap V(B'))$-rooted model of $F$ in $A'$.
    This proves that $(A',B')$ satisfies~\ref{item:find_wTS_i}-\ref{item:find_wTS_iii}.
\end{proof}

\begin{figure}[!htbp]
    \centering 
    \includegraphics[scale=1]{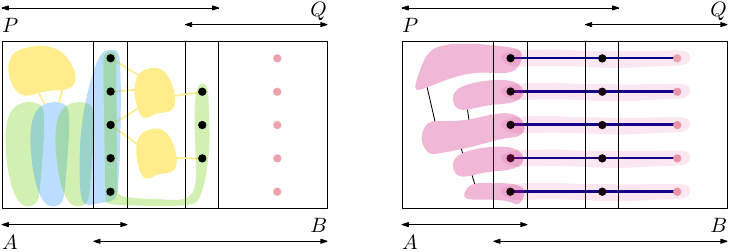} 
    \caption{An illustration of the proof of \Cref{thm:Spw}.
    On the left, we depict the situation, where $|V(P)\cap V(Q)| < |V(F)|$.
    We can extend the path decomposition by appending the bag $(V(A) \cap V(B)) \cup (V(P) \cap V(Q))$ (the last green bag in the figure).
    On the right, we depict the opposite situation, where $|V(P)\cap V(Q)| \geq |V(F)|$. Then, we simply extend the model and make it $S$-rooted.
    }
    \label{fig:excluded_path_bounded_pw}
\end{figure}

\begin{proof}[Proof of \Cref{thm:Spw}]
    The proof is illustrated in \Cref{fig:excluded_path_bounded_pw}.
    Let $F$ be a forest with at least one vertex, let $G$ be a graph, and let $S \subseteq V(G)$ such that $G$ has no $S$-rooted model of $F$.
    Suppose that $\pw(G,S) > 2|V(F)|-2$ and let $w = |V(F)|$.
    By \Cref{lemma:find_wTS_spanning_separation}, $G$ admits a separation $(A,B)$ satisfying~\ref{item:find_wTS_i}-\ref{item:find_wTS_iii}. 
    Let $(W_i \mid i\in [m])$ be a path decomposition of $(A,S\cap V(A))$ of width at most $2w-2$ with $V(A) \cap V(B) \subseteq W_m$.
    By Menger's Theorem applied to $V(A)$ and $S \cap V(B)$ there is a separation $(P,Q)$ of $G$ with $V(A) \subseteq V(P)$, $S \cap V(B) \subseteq V(Q)$ and a family $\mathcal{L}$
    of $|V(P) \cap V(Q)|$ pairwise disjoint $V(A)$--$(S \cap V(B))$ paths in $G$.
    
    First, suppose that $|V(P) \cap V(Q)| < |V(F)|$.
    Then we claim that $\big(W_1, \dots, W_m, (V(A) \cap V(B)) \cup (V(P) \cap V(Q))\big)$ is a path decomposition of $(P,S\cap V(P))$ of width at most $2w-2$
    whose last bag contains $V(P) \cap V(Q)$.
    Indeed, $|V(A) \cap V(B)| + |V(P) \cap V(Q)|-1 \leq 2w-2$.
    Let $C$ be a component of $P-\big(\bigcup_{i \in [m]} W_i \cup (V(P) \cap V(Q))\big)$.
    Since $V(A) \cap V(B) \subset W_m$, either $C$ is a component of $A-\bigcup_{i \in [m]} W_i$ or 
    $C$ is a component of $P - (V(A)\cup (V(P) \cap V(Q)))$.
    In the former case, there is $i \in [m]$ such that the neighborhood of $V(C)$ in $A$ (and so in $P$) is contained in $W_i$.
    In the latter case, the neighborhood of $V(C)$ in $P$ is contained in $(V(A) \cap V(B)) \cup (V(P) \cap V(Q))$.
    Hence, $(P,Q)$ is $(w,S)$-good, which contradicts the maximality of $(A,B)$
    because $(P,Q)$ extends $(A,B)$ and $|V(P) \cap V(Q)| < |V(F)| \leq |V(A) \cap V(B)|$.

    It follows that $|V(P) \cap V(Q)| \geq |V(F)|$.
    By \ref{item:find_wTS_ii}, there is a $(V(A) \cap V(B))$-rooted model of $F$ in $A$.
    The model combined with the paths in $\mathcal{L}$ yields an $S$-rooted model of $F$ in $G$.
    This contradicts the assumption on $G$ and ends the proof of the theorem.
\end{proof}

Finally, we proceed with a proof of \Cref{thm:lpw}.
The proof is by induction, however, since we need to keep some invariant stronger than the statement of the theorem, we encapsulate it in the following technical lemma, which we later show implies the theorem.

\begin{lemma}\label{lemma:lpw_technical}
    Let $X$ be an apex-forest with at least two vertices.
    Let $G$ be a connected graph and let $u$ be a vertex of $G$.
    If $G$ is $X$-minor-free, then $G$ has a layering $(L_j \mid j\geq 0)$ and there is a path decomposition $(W_i \mid i\in [m])$ of $G-u$ with
    \begin{enumerate}
        \item $L_0 = \{u\}$, and
        \item $|W_i \cap L_j| \leq 2|V(X)|-3$, for all $i \in [m]$ and $j \geq 1$.
    \end{enumerate}
\end{lemma}

\begin{figure}[!htbp]
    \centering 
    \includegraphics[scale=0.99]{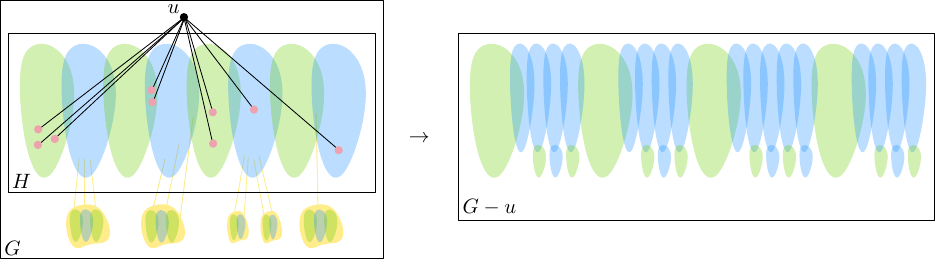} 
    \caption{An illustration of how we construct the path decomposition $(W_i \mid i \in [m])$ in the proof of \Cref{lemma:lpw_technical}.}
    \label{fig:combining_pd}
\end{figure}

\begin{proof}
    Let $x$ be a vertex of $X$ such that $X-x$ is a forest, which we denote by $F$.
    We proceed by induction on $|V(G)|$.
    If $G$ has only one vertex, then the result is clear.
    Hence, assume that $G$ has more vertices.
    
    Let $S = N(u)$ and $G' = G-u$.
    Observe that $G'$ has no $S$-rooted model of $F$, as otherwise, this model together with a branch set $\{u\}$ added would give a model of $X$ in $G$.
    By \Cref{thm:Spw}, there is a path decomposition of $(G',S)$ of width at most
    $2|V(F)|-2 = 2|V(X)|-4$.
    Let $(V_i \mid i\in [m_0])$ be such a path decomposition of $(G',S)$ with $U = \bigcup_{i \in [m_0]} V_i$ of minimum size.

    Let $C$ be a component of $G'-U$.
    We claim that $G - V(C)$ is connected. 
    Suppose to the contrary that there exists a component $C'$ of $G-V(C)$ that does not contain $u$. 
    In other words, $C'$ is disjoint from $S=N(u)$. 
    Since $G$ is connected, there is an edge $vw$ in $G$ such that $v \in V(C)$ and $w \in V(C')$.
    More precisely, $w \in U$ since otherwise $C$ is not a component of $G-U$.
    It follows that $U'=U-V(C')$ is strictly less than $U$.
    For every component $C''$ of $G'-U'$, either $C''$ is a component of $G'-U$, or $V(C'')=V(C') \cup V(C)$. 
    Since $C'$ has not neighbors in $U'$, in both cases, there exists $i \in [m_0]$ such that $N(V(C'')) \subseteq V_i \setminus V(C')$.
    Hence $(V_i \setminus V(C') \mid i \in [m_0])$ is a path decomposition of $(G',S)$.
    The width of this path decomposition is at most $2|V(X)|-4$, which contradicts the minimality of $U$.

    Let $G_C$ be obtained from $G$ by contracting $V(G) \setminus V(C)$ into a single vertex $u_C$, in particular, $G_C$ is a minor of $G$ and therefore $G_C$ is $X$-minor-free.
    Since $G$ is connected, $S$ is non-empty, thus, $|V(G_C)| \leq |V(G)| -|S \cup \{u\}|+1 \leq |V(G)|-1$.
    Hence, by induction hypothesis, there is a
    layering $(L_{C,j} \mid  j\geq 0)$ and a path decomposition $(V_{C,i} \mid i\in [m_C])$ of $G_C-u_C$ such that
    \[
     L_{C,0} = \{u_C\} \text{ and } |V_{C,i} \cap L_{C,j}| \leq 2|V(X)|-3 \text{, for every $i \in [m_C]$ and $j \geq 1$}.
    \]
    Let $L_0 = \{u\}$, $L_1 = U$, and for every $j \geq 2$, $L_j = \bigcup_C L_{C,j-1}$
    where $C$ goes over all the components of $G'-U$. See \Cref{fig:combining_layering}.
    We claim that $(L_j \mid j\geq 0)$ is a layering of $G$.
    Indeed, every edge of $G$ is either inside a layer or between two consecutive layers of $(L_j \mid j \geq 0)$ since $N(u) = S \subseteq U = L_1$, and $N(V(C)) \subseteq L_1$ and $(L_{C,j} \mid  j\geq 0)$ is a layering of $C$, for every component $C$ of $G'-U$.

\begin{figure}[!htbp]
    \centering 
    \includegraphics[scale=1]{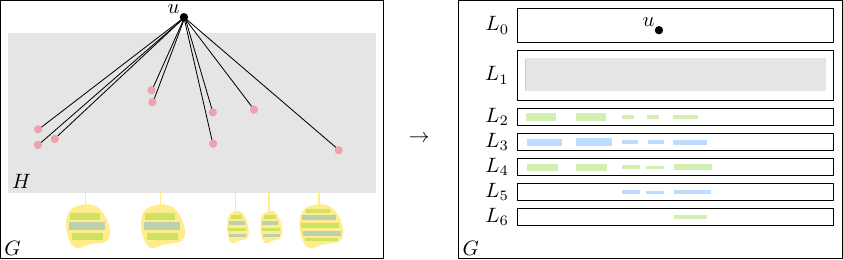} 
    \caption{An illustration of how we construct the layering $(L_j \mid j \geq 0)$ in the proof of \Cref{lemma:lpw_technical}.
    }
    \label{fig:combining_layering}
\end{figure} 

    For every component $C$ of $G'-U$, fix some $\alpha(C) \in [m_0]$ such that the neighborhood of $V(C)$ in $G$ is contained in $V_{\alpha(C)}$.
    Moreover, let the path decomposition obtained as a concatenation of the path decompositions $(V_{C,i} \mid i \in [m_C])$ for every component $C$ of $G'-U$ with $\alpha(C) = k$ be denoted by $(V_{k,i} \mid i \in [m_k])$ where $m_k = \sum_C m_C$.
    For every $k \in [m_0]$, let $V_{k,0}' = V_k$ and $V_{k,i}' = V_{k,i} \cup V_k$ for every $i\in [m_k]$.
    Observe that $(V_{k,i}' \mid 0 \leq i \leq m_k)$ is a path decomposition of the subgraph of $G'$ induced by
    $V_k \cup \bigcup_C V(C)$ where $C$ goes over all components of $G'-U$ with $\alpha(C)=k$.
    Now, let $(W_i \mid i \in [m])$ be the concatenation of the path decompositions $(V_{k,i}' \mid 0 \leq i \leq m_k)$ for each $k \in [m_0]$ in the increasing order of $k$.
    Here, $m = \sum_{k=1}^{m_0}(m_k+1)$.
    This yields a path decomposition of $G-u$.
    This construction is illustrated in \Cref{fig:combining_pd}.

    Finally, we argue that the width of $\big((W_i \mid i \in [m]), (L_j \mid j \geq 0)\big)$ is at most $2|V(X)|-3$.
    For every $i \in [m]$, we have $W_i \cap L_1  = W_i \cap U = V_k$ for some $k \in [m_0]$, and so, $|W_i \cap L_1| \leq 2|V(X)| - 3$.
    On the other hand, for every $j \geq 2$ and $i \in [m]$, we have $W_i \cap L_j = V_{C,\ell} \cap L_{C,j-1}$ for some component $C$ of $G'-U$ and $\ell \in [m_C]$, which gives $|W_i \cap L_j| \leq 2|V(X)|-3$ and ends the proof.
\end{proof}

\begin{proof}[Proof of \Cref{thm:lpw}]
    Let $X$ be an apex-forest with at least two vertices, and let $G$ be an $X$-minor-free graph.
    If $G$ has no vertex, then the result is clear.
    Hence, we assume that $V(G)$ is non-empty.
    When $G$ is connected, apply \Cref{lemma:lpw_technical} to $G$ with an arbitrary vertex $u \in V(G)$.
    We obtain a path decomposition  $(W_i \mid i \in[m])$ of $G-u$  and a layering $(L_j \mid j\geq 0)$ of $G$ such that
    $|W_i \cap L_j| \leq 2|V(X)|-3$, for every $i \in [m]$ and $j \geq 1$, and $L_0 = \{u\}$.
    Then $(W_i \cup \{u\} \mid i \in [m])$ is a path decomposition of $G$ such that every bag has intersection with every layer of $(L_j \mid j\geq 0)$ of size at most $2|V(X)|-3$.
    When $G$ is not connected, apply the above to each component of $G$ and concatenate the layerings and the path decompositions.
\end{proof}

\section{Layered treedepth}\label{sec:ltd}

In this section, we prove \Cref{thm:Std} and \Cref{thm:ltd}.
In the preliminaries section, we stated the definition of treedepth via elimination trees.
Treedepth can be also equivalently defined inductively.
Namely, treedepth of a graph is the maximum of treedepth of its components, treedepth of the one-vertex graph is $1$, and when a graph $G$ has more than one vertex and is connected, treedepth is the minimum over all vertices $v \in V(G)$ of $\td(G-v)+1$.
We claim that the treedepth of $(G,S)$ that we proposed in \Cref{sec:rooted-minor} also admits an inductive definition, which we state in terms of properties \ref{item:tdS_empty} to \ref{item:tdS_adding_a_vertex}.
Let $G$ be a graph, and $S \subset V(G)$.
\begin{enumerate}[topsep=5pt-\parskip,label={\normalfont (t\arabic*)}]
    \item If $S = \emptyset$, then $\td(G,S) = 0$. \label{item:tdS_empty}
    \item If $V(G) \neq \emptyset$, then $\td(G,S) = \max_C \td(C,S\cap V(C))$, where $C$ goes over all components of $G$. \label{item:tdS_max_connected_components}
    \item If $G$ is connected and $S \neq \emptyset$, then $\td(G,S) = 1 + \min_{u\in V(G)}\td(G-u,S-\{u\})$.  \label{item:tdS_adding_a_vertex}
\end{enumerate}
Indeed, item~\ref{item:tdS_empty} follows immediately since when $S=\emptyset$, 
we can take $F$ to be the null graph. Then, $F$ is an elimination forest of the null graph, which is a subgraph of $G$ containing $S$.
The proofs of items~\ref{item:tdS_max_connected_components} and~\ref{item:tdS_adding_a_vertex} follow by a simple induction.

Observe that we also have a monotonicity property in the following sense.
\begin{enumerate}[resume*]
    \item If $H$ is a subgraph of $G$, then $\td(H,S \cap V(H)) \leq \td(G,S)$. \label{item:tdS_subgraph}
\end{enumerate}

A \emph{depth-first-search tree}, \emph{DFS tree} for short, of $G$ is a rooted spanning tree $T$ of $G$ such that $T$ is an elimination forest of $G$.
We proceed with the proof of \Cref{thm:Std}, the key inductive step is encapsulated in the following lemma.

\begin{lemma}\label{lemma:Std_technical}
    Let $G$ be a connected graph, let $S \subseteq V(G)$, and let $T$ be a DFS tree of $G$.
    For every positive integer $\ell$, if for each root-to-leaf path $P$ in $T$,
    there are no $\ell$ pairwise disjoint $V(P)$--$S$ paths in $G$, then
    \[
      \td(G,S) \leq \binom{\ell}{2}.
    \]
\end{lemma}

\begin{figure}[!htbp]
    \centering 
    \includegraphics[scale=1]{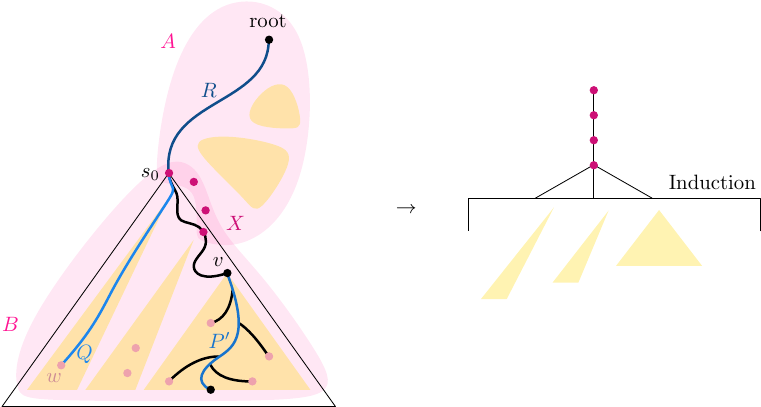} 
    \caption{
    An illustration to the proof of \Cref{lemma:Std_technical}.
    In the figure $\ell = 5$.
    On the left, we illustrate the proof by contradiction that induction can be applied to $G_v$.
    On the right, we illustrate an elimination tree that is build in the proof.
    }
    \label{fig:building_td_GS}
\end{figure}

\begin{proof}
    We proceed by induction on $\ell$.
    If $\ell=1$, then $S = \emptyset$, and so, $\td(G,S)=0$ by~\ref{item:tdS_empty}.
    Now, assume that $\ell \geq 2$.
    
    For every $u \in V(G)$, let $T_u$ be the subtree of $T$ rooted in $u$,
    and let $G_u = G[V(T_u)]$.
    Let $s_0 \in V(G)$ be the vertex with maximum depth in $T$ such that $S \subset V(T_{s_0})$.
    Let $R$ be the path from the root to $s_0$ in $T$.
    By assumption, there are no $\ell$ pairwise disjoint $V(R)$--$S$ paths.
    Hence by Menger's Theorem, there is a separation $(A,B)$ of $G$ of order at most $\ell-1$ such that $V(R) \subseteq V(A)$ and $S \subseteq V(B)$.
    In particular, every $V(R)$--$S$ path intersects $X = V(A) \cap V(B)$.

    Consider a component $C$ of $G-X$.
    If $C$ has no vertex in $S$, then $\td(C,S\cap V(C)) = 0$.
    Therefore, we assume the opposite, namely, $V(C) \subseteq V(T_{s_0})\setminus \{s_0\}$.
    It follows that there is a child $v$ of $s_0$ with $V(C) \subseteq V(T_v)$.
    The next goal is to apply induction to $G_v$ -- this step is illustrated in~\Cref{fig:building_td_GS}.
    To this end, we claim that for every root-to-leaf path $P'$ in $T_v$ there are no $\ell-1$ pairwise disjoint $V(P')$--$S$ paths in $G_v$.
    Suppose to the contrary that there is such a root-to-leaf path $P'$.
    Let $P$ be the path connecting the root of $T$ and the unique leaf in $P'$.
    By the maximality of $s_0$, $S \nsubseteq V(T_v)$, hence, there is $w \in S$ such that $w \notin V(T_v)$.
    Let $Q$ be the shortest path from $s_0$ to $w$ in $T$.
    Observe that $Q$ is a $V(R)$--$S$ path, and $Q$ disjoint from $V(T_v)$.
    Therefore, the $\ell-1$ pairwise disjoint $V(P')$--$S$ paths in $G_v$ and $Q$ form a collection of $\ell$ pairwise disjoint $V(P)$--$S$ paths in $G$, which is a contradiction.

    By inductive hypothesis applied to $G_v$ and $T_v$ we obtain $\td(G_v,S \cap V(G_v)) \leq \binom{\ell-1}{2}$.
    By repeating the above reasoning for every component of $G-X$ and \ref{item:tdS_subgraph}, this yields $\td(C,S \cap V(C)) \leq \binom{\ell-1}{2}$ for every component of $G-X$.
    In particular, by~\ref{item:tdS_max_connected_components}, $\td(G-X,S-X) \leq \binom{\ell-1}{2}$.
    Finally, by~\ref{item:tdS_adding_a_vertex},
    \begin{align*}
        \td(G,S)
        &\leq |X| + \td(G-X, S \setminus X)
        \leq (\ell-1) + \binom{\ell-1}{2} = \binom{\ell}{2}. \qedhere
    \end{align*}
\end{proof}

\begin{proof}[Proof of \Cref{thm:Std}]
    Let $\ell$ be a positive integer, let $G$ be a graph, and let $S \subseteq V(G)$.
    We can assume that $G$ is connected due to~\ref{item:tdS_max_connected_components}.
    Assume that $G$ has no $S$-rooted model of $P_\ell$, and suppose to the contrary that $\td(G,S) > \binom{\ell}{2}$.
    Then, by \Cref{lemma:Std_technical} applied with an arbitrary DFS-tree of $G$, there is a path $P$ in $G$ and $\ell$ pairwise disjoint $V(P)$--$S$ paths in $G$.
    These paths, together with $P$, give an $S$-rooted model of $P_\ell$ in $G$. 
    This is a contradiction, which ends the proof.
\end{proof}

In the second part of this section, we prove \Cref{thm:ltd}.
The proof is quite similar to the second part of the proof of \Cref{thm:lpw} in terms of structure and content.
Again, we first prove a technical lemma and then derive the theorem.

\begin{lemma}\label{lemma:ltd_technical}
    Let $X$ be a fan with at least one vertex.
    Let $G$ be a connected graph and let $u$ be a vertex of $G$.
    If $G$ is $X$-minor-free, then $G$ has a layering $(L_j \mid j\geq 0)$ and there is an elimination forest $F$ of $G-u$ with
    \begin{enumerate}
        \item $L_0 = \{u\}$, and
        \item $|V(P) \cap L_j| \leq \binom{|V(X)|-1}{2}$ for every root-to-leaf path $P$ in $F$ and for every $j \geq 1$.
    \end{enumerate}
\end{lemma}

\begin{proof}
    Let $x$ be a vertex of $X$ such that $X-x$ is a path, and let $\ell= |V(X)|-1$.
    If $\ell = 0$, then the result is vacuously true, thus, we assume that $\ell > 0$.
    We proceed by induction on $|V(G)|$.
    If $G$ has only one vertex, then the result is clear. Hence, assume that $G$ has more vertices.

    Let $S = N(u)$ and $G' = G-u$.
    Observe that $G'$ has no $S$-rooted model of $P_\ell$
    as otherwise, this model together with a branch set $\{u\}$ added would give a model of $X$ in $G$.
    By \Cref{thm:Std}, there is an elimination forest of $(G',S)$ of vertex-height at most $\binom{\ell}{2}$. 
    Let $F'$ be such a forest with $|V(F')|$ minimum.

    Let $C$ be a component of $G'-V(F')$.
    We claim that $G - V(C)$ is connected. 
    Suppose to the contrary that there exists a component $C'$ of $G-V(C)$ that does not contain $u$. 
    In other words, $C'$ is disjoint from $S=N(u)$. 
    Since $G$ is connected, there is an edge $vw$ in $G$ such that $v \in V(C)$ and $w \in V(C')$.
    More precisely, $w \in V(F')$ since otherwise, $C$ is not a component of $G-V(F')$.
    It follows that $V(F')\setminus V(C')$ is strictly less than $V(F')$.
    Let $F''$ be the forest with the vertices $V(F') = V(C')$, where for all $x,y \in V(F')$, we have $xy \in E(F'')$ whenever there is an $\{x\}$--$\{y\}$ path in $F'$ with all internal vertices in $V(C')$.
    For every component $C''$ of $G'-V(F'')$, either $C''$ is a component of $G'-V(F')$, or $V(C'')=V(C') \cup V(C)$. 
    Since $C'$ has not neighbors in $V(F'')$, in both cases, there exists a root-to-leaf path containing the neighborhood of $V(C')$ in $G'$.
    Hence $F''$ is an elimination forest of $(G',S)$.
    The vertex-height of $F'$ is at most $\binom{\ell}{2}$, which contradicts the minimality of $F'$.

    Let $G_C$ be obtained from $G$ by contracting $V(G) \setminus V(C)$ into a single vertex $u_C$, in particular, $G_C$ is a minor of $G$ and therefore $G_C$ is $X$-minor-free.
    Since $G$ is connected, $S$ is non-empty, thus, $|V(G_C)| \leq |V(G)| -|S \cup \{u\}|+1 \leq |V(G)|-1$.
    Hence, by induction hypothesis, there is a
    layering $(L_{C,j} \mid  j\geq 0)$ and an elimination forest $F_C$ of $G_C-u_C$ such that
    \[
     L_{C,0} = \{u_C\} \text{ and } |V(P) \cap L_{C,j}| \leq \binom{\ell}{2} \text{, for every root-to-leaf path $P$ in $F_C$ and $j \geq 1$}.
    \]

    Let $L_0 = \{u\}$, $L_1 = V(F')$, and for every $j \geq 2$, $L_j = \bigcup_C L_{C,j-1}$
    where $C$ goes over all components of $G'-V(F')$.
    We claim that $(L_j \mid j\geq 0)$ is a layering of $G$.
    Indeed, every edge of $G$ is either inside a layer or between two consecutive layers of $(L_j \mid j \geq 0)$ since $N(u) = S \subseteq V(F') = L_1$, and $N(V(C)) \subseteq L_1$ and $(L_{C,j} \mid  j\geq 0)$ is a layering of $C$, for every component $C$ of $G'-U$.

    Let $Z$ be the set of all leaves of $F'$, and for each $y \in Z$, let $P_y$ be the path from the root of $F'$ to $y$ in $F'$.
    For every component $C$ of $G'-V(F')$, fix some $\alpha(C) \in Z$ such that the neighborhood of $V(C)$ in $G$ is contained in $P_{\alpha(C)}$.
    Let $F$ be a forest obtained from $F'$ in the following way.
    For each component $C$ of $G-V(F')$ add $F_C$ and edges of the form $\alpha(C)x$ for every $x$ root of $F_C$.
    Let the set of roots of $F$ be the same as $F'$.
    It follows that $F$ is an elimination forest of $G-u$.

    Finally, let $P$ be a root-to-leaf path in $F$.
    We have $|V(P) \cap L_1| \leq \binom{\ell}{2}$ since $V(P) \cap L_1$ is a vertex set of a root-to-leaf path of $F'$.
    For every $j\geq 2$, $V(P) \cap L_j \subset V(F_C)$ for some component $C$ of $G'-V(F')$, which implies $|V(P) \cap L_j| = |V(P) \cap L_{C,j-1}| \leq \binom{\ell}{2}$.
    This proves the lemma.
\end{proof}

\begin{proof}[Proof of \Cref{thm:ltd}]
    Let $X$ be a fan with at least three vertices, and let $G$ be an $X$-minor-free graph.
    If $G$ has no vertex, then the result is clear.
    Hence, we assume that $V(G)$ is non-empty.
    When $G$ is connected, apply \Cref{lemma:ltd_technical} to $G$ with an arbitrary vertex $u \in V(G)$.
    We obtain an elimination forest $F$ of $G-u$  and a layering $(L_j \mid j\geq 0)$ of $G$ such that
    $|V(P) \cap L_j| \leq \binom{|V(X)-1|}{2}$, for every root-to-leaf path $P$ and for every $j \geq 1$, and $L_0 = \{u\}$.
    Let $T$ be obtained by adding $u$ to $F$ as a new root adjacent to all the roots of $F$.
    Then $T$ is an elimination tree of $G$ witnessing that $\ltd(G,S) \leq \binom{|V(X)|-1}{2}$.
    When $G$ is not connected, apply the above to each component of $G$, 
    take for $F$ the disjoint union of the elimination forests obtained for each component, and concatenate the layerings.
\end{proof}

\section{Treewidth and tangles focused on a set of vertices}\label{sec:tw}

In this section, we prove \Cref{thm:tw_and_tangle_numbers}.
Let $G$ be a graph with at least one vertex, and $S \subset V(G)$.
\Cref{lemma:upper_tw_and_tn} directly implies the upper bound, that is, $\tw(G,S) \leq 10 \max \{\tn(G,S),2\}-12$.
\Cref{lemma:lower_tw_and_tn} implies the lower bound, that is, $\tn(G,S) - 1 \leq \tw(G,S)$.

\begin{lemma}\label{lemma:upper_tw_and_tn}
    Let $k$ be an integer with $k \geq 2$.
    Let $G$ be a graph and let $S \subseteq V(G)$ be such that there is no tangle of $(G,S)$ of order $k$.
    Then for every $R \subseteq V(G)$ with $|R| \leq 7k - 8$, there is a tree decomposition $\mathcal{D}$ of $(G,S)$ of width at most $10k-12$ such that there is a bag of $\mathcal{D}$ containing $R$.
\end{lemma}

\begin{proof}
    We proceed by induction on $|V(G)|$.
    If $|V(G)| \leq 10k-11$, then the tree decomposition consisting of a single bag $V(G)$ witnesses the statement.
    Thus, we assume that $|V(G)| \geq 10k-10 \geq 7k-8$.
    By possibly adding some vertices to $R$, we assume without loss of generality that $|R| = 7k-8$.

    Let $\mathcal{T}$ be the family of all separations $(A,B)$ of $G$ of order at most $k-1$ such that $|V(A) \cap R| \leq 4k-5$.
    By assumption, $\mathcal{T}$ is not a tangle of $(G,S)$. Therefore, one of~\ref{item:T1}-\ref{item:T4} does not hold.

    If \ref{item:T1} does not hold, then there is a separation $(A,B)$ of $G$ of order at most $k-1$ such that $|V(A) \cap R| \geq 4k-4$ and $|V(B) \cap R| \geq 4k-4$.
    Then $|R| \geq |V(A) \cap R| + |V(B) \cap R| - |V(A) \cap V(B)| \geq 8k-8 - (k-1) = 7k-7 > 7k-8$, a contradiction.

    If \ref{item:T2} does not hold, then there are separations $(A_1,B_1),(A_2,B_2),(A_3,B_3)$ in $\mathcal{T}$ such that $A_1 \cup A_2 \cup A_3 = G$.
    Let $Z= \bigcup_{i=1}^3 \big(V(A_i) \cap V(B_i)\big)$. 
    Let $C$ be a component of $G-Z$, let $G_C = G[V(C) \cup N(C)]$, and let $R_C = N(V(C)) \cup (R \cap V(C))$.
    Since $V(C) \subseteq V(A_i)$ for some $i \in \{1,2,3\}$, $|V(C) \cap R| \leq |V(A_i) \cap R| \leq 4k-5$.
    Observe that $V(G_C) = V(C) \cup N(C) \subset A_i$, and thus, 
    \[|V(G) \setminus V(G_C)| \geq |V(B_i) \setminus V(A_i)| \geq |(V(B_i) \setminus V(A_i)) \cap R| \geq 3k-3 > 0.\]
    Moreover, since $N(V(C)) \subset Z$, $|N(V(C))| \leq |Z| = 3(k-1)$.
    Hence, $|R_C| \leq |V(C) \cap R| + |N(V(C))| \leq 4k-5+3k-3 = 7k - 8$.
    In order to apply induction to $G_C$ and $R_C$, we have to argue that $|V(G_C)| < |V(G)|$.
    By induction hypothesis applied to $G_C$ and $R_C$, there is a tree decomposition $\big(T_C,(W_{C,x} \mid x \in V(T_C))\big)$ of $(G_C,S\cap V(G_C))$ of width at most $10k-12$ such that $R^C \subseteq W_{C,r_C}$ for some $r_C \in V(T_C)$.
    Let $T$ be obtained from the disjoint union of $T_C$ for all components $C$ of $G-Z$ by adding a new vertex $r$ and edges $rr_C$ for every component $C$ of $G-Z$.
    Finally, let $W_r = Z \cup R$, and for every component $C$ of $G-Z$ and every $x \in V(T_C)$, let $W_x = W_{C,x}$.
    Observe that $|W_r| \leq |Z|+|R| \leq 3(k-1)+7k-8 = 10k-11$.
    Every component of $G - \bigcup_{x \in V(T)} W_x$ is a subgraph of a component of $G_C - \bigcup_{x \in V(T_C)} W_{C,x}$ for some component $C$ of $G-Z$.
    Therefore, $\big(T,(W_x \mid x \in V(T))\big)$ is a tree decomposition of $(G,S)$ of width at most $10k-12$ such that $R \subseteq W_r$.

    If \ref{item:T3} does not hold, then there is a separation $(A,B) \in \mathcal{T}$ such that $V(A)=V(G)$.
    It follows that $|R| = |R \cap V(A)| \leq 4k-5 < 7k-8 = |R|$, a contradiction.
    
    If \ref{item:T4} does not hold, then there is a separation $(A,B) \in \mathcal{T}$ such that $S \subseteq V(A)$.
    Let $R' = (R \cap V(A)) \cup (V(A) \cap V(B))$.
    Observe that $|R'| \leq 4k-5 + (k-1) =5k-6 \leq 7k-8$.
    By induction hypothesis applied to $A$ and $R'$, there is tree decomposition $\big(T',(W_x \mid x \in V(T'))\big)$ of $(A,S\cap V(A))$ of width at most $10k-12$ such that $R' \subseteq W_{r'}$ for some $r' \in V(T')$.
    Let $T$ be obtained from $T'$ by adding a new vertex $r$ and the edge $rr'$.
    Finally, set $W_r = R \cup (V(A) \cap V(B))$ and observe that $|W_r|\leq |R|+k-1\leq 7k-8+k-1 = 8k-9 \leq 10k-12$.
    Every component of $G - \bigcup_{x \in V(T)} W_x$ is either a component of $B-A$ or is a subgraph of a component of $A - \bigcup_{x \in V(T')} W_{x}$.
    In both cases, the neighborhood of the component is contained in a single bag.
    Therefore, $\big(T,(W_x \mid x \in V(T))\big)$ is a tree decomposition of $(G,S)$ of width at most $10k-12$ such that $R \subseteq W_r$.
\end{proof}

In the proof of \Cref{lemma:lower_tw_and_tn} we use the following simple observation.

\begin{obs}\label{obs:tangles}
    Let $k$ be a positive integer, let $G$ be a graph, and let $\mathcal{T}$ be a tangle of $G$ of order $k$.
    Let $(A,B)$ and $(A',B')$ be two separations of $G$ of order at most $k-1$ such that $V(A) = V(A')$ and $V(B) = V(B')$.
    Then, 
    \[ (B,A) \notin \mathcal{T} \Longleftrightarrow (A,B) \in \mathcal{T} \Longleftrightarrow (A',B') \in \mathcal{T}  \Longleftrightarrow (B',A') \notin \mathcal{T}.\]
\end{obs}
\begin{proof}
    The first and last equivalences are clear by \ref{item:T1} and \ref{item:T2}.
    In order to prove the middle equivalence, suppose to the contrary that $(A,B) \in \mathcal{T}$ and $(A',B') \not\in \mathcal{T}$.
    By~\ref{item:T2}, $(B',A') \in \mathcal{T}$.
    Observe that $(A \cup B', A' \cap B)$ is a separation of $G$, and its order is at most $k-1$.
    By~\ref{item:T3}, $(A \cup B', A' \cap B) \not\in \mathcal{T}$, hence by~\ref{item:T1},  $(A' \cap B, A \cup B') \in \mathcal{T}$.
    But then $B' \cup A \cup (A' \cap B) = G$, which contradicts~\ref{item:T3}.    
\end{proof}

\begin{lemma}\label{lemma:lower_tw_and_tn}
    Let $k$ be a positive integer, let $G$ be a graph with at least one vertex, and let $S \subseteq V(G)$.
    If $\tn(G,S) \geq k$, then $\tw(G,S) \geq k-1$.
\end{lemma}

\begin{proof}
    Let $\mathcal{T}$ be a tangle of $(G,S)$ of order $k$.
    Suppose to the contrary that there is a tree decomposition $\big(T_0,(W_x \mid x \in V(T_0))\big)$ of $(G,S)$ of width at most $k-2$. Let $U=\bigcup_{x\in V(T_0)} W_x$. 
    By possibly adding some vertices to $T_0$ and some bags to $(W_x \mid x \in V(T_0))$,
    without loss of generality we can assume that every vertex in $U$ is in at least two bags.
    For every component $C$ of $G-U$, there is a bag $x_C \in V(T_0)$ such that $N(V(C)) \subset W_{x_C}$.
    Let $T$ be obtained from $T_0$ by adding a new vertex $u_C$ and the edge $x_Cu_C$ for every component $C$ of $G-U$.
    Let $W_{u_C} = N(V(C)) \cup V(C)$ for every component $C$ of $G-U$.
    It follows that $\big(T,(W_x \mid x \in V(T))\big)$ is a tree decomposition of $G$.
    While this tree decomposition may have large width, for every edge $xy \in E(T)$, we have $|W_x \cap W_y| \leq k-1$.
    
    Let $Z = V(T) \setminus V(T_0)$ be the set of all added vertices.
    For every $uv \in E(T)$, let $T_{u\mid v}$ be the component of $T \setminus uv$ containing $u$, and let $G_{u\mid v}$ be the subgraph $G[\bigcup_{x \in V(T_{u\mid v})} W_x]$.

    Let $\vec{T}$ be the directed graph with the vertex set $V(T)$ and the arc set consisting of all the pairs $(u,v) \in V(T)^2$ such that $uv \in E(T)$ and for every separation $(A,B)$ of $G$ with $V(A) = V(G_{u \mid v})$ and $V(B) = V(G_{v\mid u})$, we have $(A,B) \in \mathcal{T}$.
    By \Cref{obs:tangles}, $\vec{T}$ is an orientation of $T$.
    Since $T$ is a tree, $\vec{T}$ is acyclic, and thus, there is a sink $x$ in $\vec{T}$.
    If $x \in Z$, then the neighbor $y$ of $x$ in $T$ is such that $(G_{y\mid x},G_{x \mid y} \setminus E(G_{y \mid x})) \in \mathcal{T}$, which contradicts~\ref{item:T4} since $S \subseteq V(G_{y \mid x})$.
    Hence, $x \not\in Z$, and so, $|W_x| \leq k-1$.
    Let $y_1, \dots, y_d$ be the neighbors of $x$ in $T$.
    For every $i \in [d]$, let $(A_i, B_i)$ be a separation of $G$ with $A_i = G\left[\bigcup_{j=1}^i V(G_{y_j \mid x})\right]$ and $B_i = G\left[\bigcap_{j=1}^i V(G_{x \mid y_j})\right] - E(A_i)$.
    It follows that $V(A_i) \cap V(B_i) \subset W_x$.
    Therefore, $(A_i,B_i)$ has order at most $|W_x|\leq k-1$ for every $i \in [d]$.
    
    We claim that $(A_i,B_i) \in \mathcal{T}$ for every $i \in [d]$.
    We prove this by induction on $i$.
    The fact that $x$ is a sink implies that $(A_1,B_1) \in \mathcal{T}$. 
    For the inductive step, let $1 < i \leq d$, and assume that $(A_{i-1},B_{i-1}) \in \mathcal{T}$.
    Suppose to the contrary that $(B_i,A_i) \in \mathcal{T}$.
    Since $(A_{i-1},B_{i-1}) \in \mathcal{T}$ and $(G_{y_i \mid x},G_{x \mid y_i} - E(G_{y_i \mid x})) \in \mathcal{T}$ by~\ref{item:T2}, $A_{i-1} \cup G_{y_i \mid x} \cup B_i \neq G$, which is false since $A_i = A_{i-1} \cup G_{y_i \mid x}$, and yields $(A_i,B_i) \in \mathcal{T}$.

    The above in particular, implies that $(A_d,B_d) \in \mathcal{T}$.
    However, since every vertex in $U$ is in at least two bags, 
    we have $V(G) = \bigcup_{z \in V(T) \setminus \{x\}} W_z = V(A_d)$,
    which contradicts~\ref{item:T3}, and shows that there is no tree decomposition of $(G,S)$ of width at most $k-2$.
\end{proof}

\section{A lower bound for \texorpdfstring{\Cref{cor:td-diam,cor:pw-diam}}{Corollaries~\ref{cor:td-diam} and~\ref{cor:pw-diam}}}\label{sec:lower-bounds}

Let $G$ be a graph.
Recall that the radius of a graph is the minimum over all vertices $u\in V(G)$ of $\max_{v\in V(G)}\dist_G(u,v)$ and that it is at least half of the diameter of $G$.

\begin{theorem}\label{thm:td_rad_lowerbound}
    Let $\ell$ and $r$ be integers such that $\ell\geq 2$ and $r\geq 0$. There is a fan $X$ on at least $\ell+1$ vertices and an $X$-minor-free graph $G$ with radius at least $r$ such that 
    \[
    \td(G) - 1 \geq \pw(G) \geq \left\lfloor\frac{\ell}2 \right\rfloor \left(r- \left\lfloor\frac{\ell}{2} \right\rfloor \right).
    \]
\end{theorem}

Note that~\cref{cor:pw-diam} gives a matching upper bound $\pw(G)=\Oh(\ell r)$ and
\cref{cor:td-diam} gives an upper bound $\td(G)=\Oh(\ell^2 r)$.

\begin{figure}[!htbp]
    \centering 
    \includegraphics[scale=0.94]{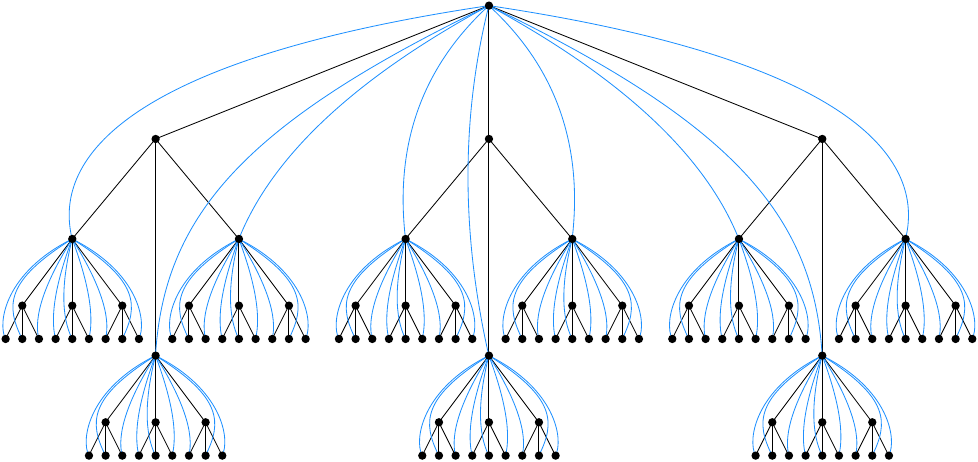} 
    \caption{The construction of $G$ in the proof of \Cref{thm:td_rad_lowerbound}, for $\ell = k = 4$ and $r = 2$.
    } \label{fig:lower_bound}
\end{figure}

\begin{proof}
    Let $X$ be the graph obtained from the path on $\ell$ vertices by adding a universal vertex $x$.
    Let $k$ be an integer such that $k \geq \frac12 \ell r$, and let $T$ be a rooted complete ternary tree of vertex-height $k+1$. 
    Let $G=T$ when $\ell < 4$.
    When $\ell \geq 4$, let $G$ be obtained from $T$ in the following way.
    For every integer $i$ such that $0\leq i\leq (1 \slash \left\lfloor \frac{\ell}{2} \right\rfloor) k-1$, for every vertex $u$ at depth $1+\left\lfloor \frac{\ell}{2} \right\rfloor i$ and every vertex $v$ at depth $1+\left\lfloor \frac{\ell}{2} \right\rfloor (i+1)$ such that $u$ is an ancestor of $v$, we add the edge $uv$.
    See \Cref{fig:lower_bound}.

    First, we show that $G$ is $X$-minor-free.
    If $\ell < 4$, then $X$ contains a cycle and $G$ is a tree, which cannot have $X$ as a minor.
    Thus, we assume that $\ell\geq 4$. 
    Suppose to the contrary that $G$ contains a model of $X$, and let $\mathcal{X}$ be such a model that is inclusion-wise minimal.
    Since $X$ is $2$-connected, $\mathcal{X}$ must be contained in a $2$-connected subgraph of $G$.
    Consider a maximal $2$-connected subgraph $H$ of $G$ containing $\mathcal{X}$.
    Let $u$ be the vertex in $V(H)$ with minimum depth in $T$. 
    By construction of $G$, $H-u$ is a rooted complete ternary tree $T'$ with vertex-height $\left\lfloor\frac{\ell}{2} \right\rfloor$, the root of $T'$ is the only child of $u$ in $V(H)$, and $uv$ is an edge of $H$ for every leaf $v$ of $T'$. 
    The branch set of $x$ in $\mathcal{X}$ must contain $u$, as otherwise, $H-u = T'$ does not contain a model of $K_3$. 
    Therefore, $T'$ must contain a model of a path on $\ell$ vertices, but the longest path in $T'$ has only $2\left\lfloor \frac{\ell}{2} \right\rfloor-1 < \ell$ vertices. 
    This is a contradiction, thus, $G$ is $X$-minor-free.

    By symmetry of the construction, the radius of $G$ is witnessed by the root of $G$.
    Let $t$ be the root of $G$, and let $v$ be a vertex in $T$ of depth $i$ for some integer $i$ with $1 \leq i \leq k+1$.
    We have, $\dist_G(t,v) \leq (1 \slash \left\lfloor \frac{\ell}{2} \right\rfloor) (i-1) + \left\lfloor \frac{\ell}{2} \right\rfloor \leq (1 \slash \left\lfloor \frac{\ell}{2} \right\rfloor) k + \left\lfloor \frac{\ell}{2} \right\rfloor$, hence, the radius of $G$ is at most $(1 \slash \left\lfloor \frac{\ell}{2} \right\rfloor) k + \left\lfloor \frac{\ell}{2} \right\rfloor$.
    On the other hand, when $v$ is a leaf, $\dist_G(t,v)\geq (1 \slash \left\lfloor \frac{\ell}{2} \right\rfloor) k \geq r$ by definition of $k$.
    So, the radius of $G$ is at least $r$.
    We obtain $r \leq (1 \slash \left\lfloor \frac{\ell}{2} \right\rfloor) k + \left\lfloor \frac{\ell}{2} \right\rfloor$.
    
    Since $G$ contains a rooted complete ternary tree of vertex-height $k+1$ as a subgraph, $\pw(G)\geq k$. Moreover, $\td(G) - 1 \geq \pw(G)$.
    Therefore $G$ is an $X$-minor-free graph $G$ with radius of at least $r$ witnessing the assertion of the theorem. \qedhere
\end{proof}

\section{Erd\H{o}s-P\'osa property}\label{sec:EP-property}

In this section, we discuss the applications of our techniques to Erd\H{o}s-P\'osa properties for rooted models.
We start with a classical statement by Robertson and Seymour on families of connected subgraphs in graphs of bounded treewidth.

\begin{lemma}[\cite{GM5}, Statement (8.7)]\label{lemma:RS-EP-property}  
Let $G$ be a graph, let $\mathcal{W}=(T,(W_x \mid x\in V(T)))$ be a tree decomposition of $G$, and let $\mathcal{F}$ be a family of connected subgraphs of $G$.
For every positive integer $k$, either
    \begin{enumerate}
        \item there are $k$ pairwise vertex-disjoint subgraphs in $\mathcal{F}$ or 
        \item there is a set $Z \subseteq V(G)$ that is the union of at most $k-1$ bags of $\mathcal{W}$ such that $V(F) \cap Z \neq \emptyset$ for every $F \in \mathcal{F}$.
    \end{enumerate}
\end{lemma}

It turns out that the analog version for treewidth focused on a prescribed set of vertices holds.

\begin{lemma}\label{lemma:EP_Stw}  
Let $G$ be a graph, let $S\subseteq V(G)$, let $\mathcal{W}=(T,(W_x \mid x\in V(T)))$ be a tree decomposition of $(G,S)$, and let $\mathcal{F}$ be a family of connected subgraphs of $G$ each of them intersecting $S$.
For every positive integer $k$, either
    \begin{enumerate}
        \item there are $k$ pairwise vertex-disjoint subgraphs in $\mathcal{F}$ or \label{item:EP_i}
        \item there is a set $Z \subseteq V(G)$ that is the union of at most $k-1$ bags of $\mathcal{W}$ such that $V(F) \cap Z \neq \emptyset$ for every $F \in \mathcal{F}$. \label{item:EP_ii}
    \end{enumerate}
\end{lemma}

\begin{proof}
    Let $k$ be a positive integer, and suppose that~\ref{item:EP_i} does not hold.
    For every $F \in \mathcal{F}$, let $T_F = T[\{x \in V(T) \mid W_x \cap V(F) \neq \emptyset\}]$.
    It follows that for every $F \in \mathcal{F}$, $T_F$ is a non-empty subtree of $T$.

    We claim that if $F_1,F_2 \in \mathcal{F}$ and $V(F_1) \cap V(F_2) \neq \emptyset$, then $V(T_{F_1}) \cap V(T_{F_2}) \neq \emptyset$.
    Indeed, if $u \in V(F_1) \cap V(F_2)$,
    then either $u \in W_x$ for some $x \in V(T)$ and so $x \in V(T_{F_1}) \cap V(T_{F_2})$, or $u \in V(C)$ for some component $C$ of $G- \bigcup_{x \in V(T)} W_x$.
    Then, since $\mathcal{W}$ is a tree decomposition of $(G,S)$, there exists $x \in V(T)$
    such that $N(V(C)) \subseteq W_x$. Moreover, $V(F_i) \cap S \neq \emptyset$, and so, $V(F_i) \cap N(V(C)) \neq \emptyset$, for each $i \in \{1,2\}$. Hence $x \in V(T_{F_1}) \cap V(T_{F_2})$.

    Since~\ref{item:EP_i} is false, we deduce that there are no $k$ disjoint members of $\{T_F \mid F \in \mathcal{F}\}$. Then by Helly property for subtrees of $T$, there are $k-1$ bags of $\cgW$ whose union $Z$ intersects every member of $\mathcal{F}$. Therefore, \ref{item:EP_ii} holds.
\end{proof}

\Cref{lemma:EP_Stw} with \Cref{thm:Stw} yield that outer-rooted models of a fixed connected plane graph admit the Erd\H{o}s-P\'osa property.
Recall that $\gm$ is the minimum function such that for every positive integer $\ell$, if a graph $G$ has no model of $\boxplus_\ell$, then $\tw(G) \leq \gm(\ell)$.

\begin{cor}
    For every connected plane graph $H$, 
    for every graph $G$, 
    for every $S\subseteq V(G)$, and
    for every positive integer $k$, 
    either
    \begin{enumerate}
        \item $G$ has $k$ vertex-disjoint $S$-outer-rooted models of $H$ or\label{item:EP_thm_i}
        \item there exists a set $Z \subseteq G$ such that $|Z| \leq 3 (k-1)(\gm(98304 \cdot k^4|V(H)|^4)+2)$ and $G-Z$ has no $S$-outer-rooted model of $H$.\label{item:EP_thm_ii}
    \end{enumerate}
\end{cor}

\begin{proof}
    Let $H$ be a connected plane graph.
    For every positive integer $k$, let $k\cdot H$ denote the plane graph consisting of $k$ disjoint copies of $H$ drawn in the plane in such a way that the outer face of each copy belongs to the outer face.
    Suppose that $G$, $S$, $k$ verify the conditions in the statement. 
    Assume that~\ref{item:EP_thm_i} does not hold.
    In other words, $G$ has no $S$-outer-rooted model of $k \cdot H$.
    Therefore, $\tw(G,S) \leq 3 \gm(98304\cdot k^4|V(H)|^4)+1$ by \Cref{thm:Stw}.
    Then by \Cref{lemma:EP_Stw} applied to the family of all the connected subgraphs of $G$ containing an $S$-outer-rooted model of $H$, there exists a set $Z$ of at most $(\tw(G,S)+1)(k-1)$ vertices in $G$ such that $G-Z$ has no $S$-outer-rooted-model of $H$.
\end{proof}

Recently, Dujmović, Joret, Micek, and Morin~\cite{dujmović2024tight} showed that for every tree $T$, for every graph $G$, for every positive integer $k$, either $G$ has $k$ disjoint models of $T$, or there is a set $Z$ of at most $|V(T)|(k-1)$ vertices such that $G-Z$ is $T$-minor-free.
\Cref{thm:Spw} and \Cref{lemma:EP_Stw} imply the following Erd\H{o}s-P\'osa property for rooted models of trees.

\begin{cor}\label{thm:EP-for-rooted-trees}
    For every tree $T$, 
    for every graph $G$, 
    for every $S\subseteq V(G)$, and
    for every positive integer $k$,  
    either
    \begin{enumerate}
        \item $G$ has $k$ vertex-disjoint $S$-rooted models of $T$ or
        \item there exists a set $Z \subseteq G$ such that $|Z| \leq (2k|V(T)|-1) (k-1)$ and $G-Z$ has no $S$-rooted model of $T$.
    \end{enumerate}
\end{cor}

\section{Open problems}\label{sec:open}

Some of the bounds that we provided are not tight, we summarize potential improvements below.

\begin{problem}
Within~\Cref{thm:Std}, 
we show that
    for every positive integer $\ell$, for every graph $G$, and for every $S\subseteq V(G)$, if $G$ has no $S$-rooted model of $P_\ell$, then $\td(G,S)\leq \binom{\ell}{2}$.
    Is there a better bound? Perhaps linear in $\ell$?
\end{problem}

\begin{problem}
    Within \Cref{thm:Spw} we show that for every forest $F$ with at least one vertex,
    for every graph $G$, for every $S \subseteq V(G)$, if $G$ has no $S$-rooted model of $F$, then $\pw(G,S) \leq 2|V(F)|-2$.
    Is there a better bound? Perhaps $|V(F)|-2$?
\end{problem}

Any improvement of the bounds in \Cref{thm:Std,thm:Spw} yields improvement of the bounds in \Cref{thm:ltd,thm:lpw} respectively.
The lower bounds in \Cref{thm:lpw,thm:ltd} are witnessed by taking a clique on $|V(X)|-1$ vertices: such a graph is $X$-minor-free regardless of the structure of $X$ and $\lpw(K_{|V(X)|-1}) = \ltd(K_{|V(X)|-1}) = \lceil (|V(X)|-1)\slash 2\rceil$.

\begin{problem}
Within~\Cref{thm:EP-for-rooted-trees} we show the Erd\H{o}s-P\'osa property for $S$-rooted trees with a bound $(2k|V(T)|-1)(k-1)=\mathcal{O}(k^2)|V(T)|$. 
Is there a better bound? 
Perhaps $\mathcal{O}(k)|V(T)|$?
\end{problem}

\subsection*{Acknowledgments}
We are grateful to two anonymous reviewers for their valuable feedback on a previous version of this article.

\bibliographystyle{plain}
\bibliography{biblio}

@article{Marx2017,
  title = {Rooted grid minors},
  volume = {122},
  ISSN = {0095-8956},
  url = {http://dx.doi.org/10.1016/j.jctb.2016.07.003},
  DOI = {10.1016/j.jctb.2016.07.003},
  journal = {Journal of Combinatorial Theory,  Series B},
  publisher = {Elsevier BV},
  author = {Marx,  Daniel and Seymour,  Paul and Wollan,  Paul},
  year = {2017},
  pages = {428–437},
  note = {\href{https://arxiv.org/abs/1307.8138}{arXiv:1307.8138}}
}

@article{seymour2023shorter,
  title={A shorter proof of the path-width theorem}, 
  author={Paul Seymour},
  year={2023},
  journal={{arXiv preprint}},
  note = {\href{https://arxiv.org/abs/2309.05100}{arXiv:2309.05100}}
}

@article{Diestel1995,
  title = {{Graph Minors I: A Short Proof of the Path-width Theorem}},
  volume = {4},
  ISSN = {1469-2163},
  url = {http://dx.doi.org/10.1017/S0963548300001450},
  DOI = {10.1017/s0963548300001450},
  number = {1},
  journal = {Combinatorics,  Probability and Computing},
  publisher = {Cambridge University Press (CUP)},
  author = {Diestel,  Reinhard},
  year = {1995},
  pages = {27–30}
}

@ARTICLE{GM1,
  author =	 {Robertson, Neil and Seymour, Paul},
  title =	 {Graph minors. {I}. {E}xcluding a forest},
  journal =	 {Journal of Combinatorial Theory, Series B},
  year =	 1983,
  volume =	 35,
  pages =	 {39--61},
  number =	 1,
MRNUMBER = {723569},
  DOI =		 {10.1016/0095-8956(83)90079-5},
}

@article{GM5,
title = {Graph minors. {V}. {Excluding a planar graph}},
journal = {Journal of Combinatorial Theory, Series B},
volume = {41},
number = {1},
pages = {92-114},
year = {1986},
issn = {0095-8956},
doi = {https://doi.org/10.1016/0095-8956(86)90030-4},
author = {Neil Robertson and Paul Seymour},
}

@book{Sparsity-book,
  AUTHOR =	 {Ne{\v{s}}et{\v{r}}il, Jaroslav and Ossona de Mendez,
                  Patrice},
  TITLE =	 {Sparsity},
  SERIES =	 {Algorithms and Combinatorics},
  VOLUME =	 28,
  PUBLISHER =	 {Springer},
  YEAR =	 2012,
  MRNUMBER =	 2920058,
  URL =  {http://dx.doi.org/10.1007/978-3-642-27875-4},
}

@article{layered-treewidth,
title = {Layered separators in minor-closed graph classes with applications},
journal = {Journal of Combinatorial Theory, Series B},
volume = {127},
pages = {111-147},
year = {2017},
issn = {0095-8956},
doi = {https://doi.org/10.1016/j.jctb.2017.05.006},
url = {https://www.sciencedirect.com/science/article/pii/S0095895617300473},
author = {Vida Dujmović and Pat Morin and David R. Wood},
note = {\href{https://arxiv.org/abs/1306.1595}{arXiv:1306.1595}}
}

@article{layered-pathwidth,
author = {Dujmovi\'{c}, Vida and Eppstein, David and Joret, Gwena\"{e}l and Morin, Pat and Wood, David R.},
title = {Minor-Closed Graph Classes with Bounded Layered Pathwidth},
journal = {SIAM Journal on Discrete Mathematics},
volume = {34},
number = {3},
pages = {1693-1709},
year = {2020},
doi = {10.1137/18M122162X},
URL = {https://doi.org/10.1137/18M122162X},
eprint = {https://doi.org/10.1137/18M122162X},
note = {\href{https://arxiv.org/abs/1810.08314}{arXiv:1810.08314}}
}

@article{quickly-excluding-a-forest,
title = {Quickly excluding a forest},
journal = {Journal of Combinatorial Theory, Series B},
volume = {52},
number = {2},
pages = {274-283},
year = {1991},
issn = {0095-8956},
doi = {https://doi.org/10.1016/0095-8956(91)90068-U},
author = {Dan Bienstock and Neil Robertson and Paul Seymour and Robin Thomas},
}

@article{paul2023universal,
      title={Universal Obstructions of Graph Parameters}, 
      author={Christophe Paul and Evangelos Protopapas and Dimitrios M. Thilikos},
      year={2023},
      journal={arXiv preprint},
      note={\href{https://arxiv.org/abs/2304.14121}{arXiv:2304.14121}}
}

@article{dujmović2024tight,
      title={{Tight bound for the Erd\H{o}s-P\'osa property of tree minors}}, 
      author={Vida Dujmović and Gwenaël Joret and Piotr Micek and Pat Morin},
      year={2024},
      journal = {arXiv preprint},
      note = {\href{https://arxiv.org/abs/2403.06370}{arXiv:2403.06370}}
}

@article{Bose2022,
  title = {Separating layered treewidth and row treewidth},
  volume = {24(1)},
  ISSN = {1365-8050},
  url = {http://dx.doi.org/10.46298/dmtcs.7458},
  DOI = {10.46298/dmtcs.7458},
  journal = {Discrete Mathematics \& Theoretical Computer Science},
  publisher = {Centre pour la Communication Scientifique Directe (CCSD)},
  author = {Bose,  Prosenjit and Dujmović,  Vida and Javarsineh,  Mehrnoosh and Morin,  Pat and Wood,  David R.},
  year = {2022},
  note = {\href{https://arxiv.org/abs/2105.01230v3}{arXiv:2105.01230v3}}
}

@article{Robertson1994,
  title={Quickly excluding a planar graph},
  author={Robertson, Neil and Seymour, Paul and Thomas, Robin},
  journal={Journal of Combinatorial Theory, Series B},
  volume={62},
  number={2},
  pages={323--348},
  year={1994},
  publisher={Elsevier}
}

@misc{kawarabayashi2021quicklyexcludingnonplanargraph,
      title={Quickly excluding a non-planar graph}, 
      author={Ken-ichi Kawarabayashi and Robin Thomas and Paul Wollan},
      year={2021},
      eprint={2010.12397},
      archivePrefix={arXiv},
      primaryClass={math.CO},
      url={https://arxiv.org/abs/2010.12397}, 
  note = {\href{https://arxiv.org/abs/2010.12397}{arXiv:2010.12397}}
}

@article{Chuzhoy2021,
  title = {Towards tight(er) bounds for the Excluded Grid Theorem},
  volume = {146},
  ISSN = {0095-8956},
  url = {http://dx.doi.org/10.1016/j.jctb.2020.09.010},
  DOI = {10.1016/j.jctb.2020.09.010},
  journal = {Journal of Combinatorial Theory,  Series B},
  publisher = {Elsevier BV},
  author = {Chuzhoy,  Julia and Tan,  Zihan},
  year = {2021},
  month = jan,
  pages = {219–265}
}

\appendix

\section{Outer-rooted models of plane graphs}
\label{appendix:models}

In this appendix, we show that 
for every plane graph $H$, for every graph $G$ and $S \subset V(G)$, 
if $G$ contains an $S$-outer-rooted model of $\boxplus_{2|V(H)|}$, 
then $G$ contains an $S$-outer-rooted model of~$H$.
This shows that~\Cref{thm:Stw-grids} implies~\Cref{thm:Stw}.

Let $k$ and $\ell$ be positive integers.
Consider the $k\times \ell$ grid.
For each $i\in [k]$, we call the subset $\{(i,1),(i,2),\dots,(i,\ell)\}$ of vertices of the grid a \emph{row}.
For each $j \in [\ell]$, we call the subset $\{(1,j),(2,j),\dots,(k,j)\}$ a \emph{column}. 
The rows and columns are naturally ordered.

The following result follows from the proof of~\cite[statement~1.3]{Robertson1994}.

\begin{lemma}\label{lemma:RWT_planar_hamiltonian}
    Let $n$ be a positive integer and let $H$ be a Hamiltonian $n$-vertex planar graph.
    For every $u \in V(H)$, there is a model of $H$ in $\boxplus_n$
    such that the branch set of $u$ is contained in the union of the first row and the first column of $\boxplus_n$.
\end{lemma}

\begin{proof}[Proof of \Cref{thm:Stw} assuming~\Cref{thm:Stw-grids}]
    Let $H$ be a plane graph.
    By possibly adding some edges without changing the vertex set of the outer face, we assume that $H$ is connected.
    Let $H'$ be a plane graph obtained by adding to $H$ a vertex $u$ adjacent to all the vertices of the outer face and placing it in the outer face.
    Since $H$ is connected,
    there is a spanning tree $T$ of $H'$ in which $u$ is a leaf.
    Now, the goal is to transform $H'$ into a Hamiltonian plane graph $H''$ containing $H$ as a minor.
    We first describe the construction informally and later give it in full detail.
    We replace every non-leaf vertex of $T$ by a cycle of length equals to that vertex degree and duplicate each edge of $T$. We draw these new vertices and edges along the original drawing of $T$ to keep planarity.
    In other words, we "cut open" $T$.
    The duplicated edges of $T$ in $H''$ form a Hamiltonian cycle of the graph.
    See \Cref{fig:hamiltonian}.

    \begin{figure}[!htbp]
    \begin{center}
        \includegraphics{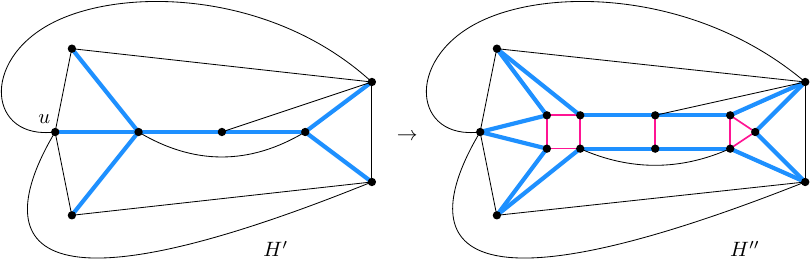}
    \end{center}
    \caption{We depict the construction in the proof of~\Cref{lemma:RWT_planar_hamiltonian}.
    The vertex $u$ is adjacent to all the vertices of the outer face of $H = H'-u$.
    We mark in blue a spanning tree $T$ of $H'$ such that $u$ is a leaf.
    On the right-hand-side, we obtain $H''$ by duplicating the blue edges and adding a pink cycle for each original non-leaf vertex in $T$.}
    \label{fig:hamiltonian}
    \end{figure}
    
    Next, let us proceed with the formal description of the construction.
    For every $v \in V(H')$, let $\pi_v$ be the cyclic permutation 
    of $N_{H'}(v)$ given by the clockwise order of $N_{H'}(v)$ around $v$.
    For every $v \in V(H')$ and every $w \in N_{H'}(v)$, 
    we define $\pi^T_v(w)$ to be the first vertex $w' \in N_T(v)$ after $w$ along $\pi_v$ such that $vw' \in E(T)$.
    The cyclic order of the neighbors of $v$ in $T$ given by the plane embedding of $T$ inherited from $H'$ is $\pi^T_v\vert_{N_T(v)}$.
    Let $H''$ be the graph defined by
    \begin{align*}
        V(H'') &=
        \{(v,w) \mid v\in V(T), w \in N_{T}(v)\} \\
    \intertext{and}
        E(H'') &= \{(v,w)(v,\pi^T_v(w)) \mid v \in V(T), w \in N_T(v), w \neq \pi^T_v(w)\} \\
        & \hspace{5.5mm} \cup \{(v,w)(w,\pi_w^T(v)),\; (v,\pi_v^T(w))(w,v) \mid vw \in E(T)\} \\
        & \hspace{5.5mm} \cup \{(v,\pi_v^T(w))(w,\pi^T_w(v)) \mid vw \in E(H') \setminus E(T)\}.
    \end{align*}

    \begin{figure}[!htbp]
    \begin{center}
        \includegraphics{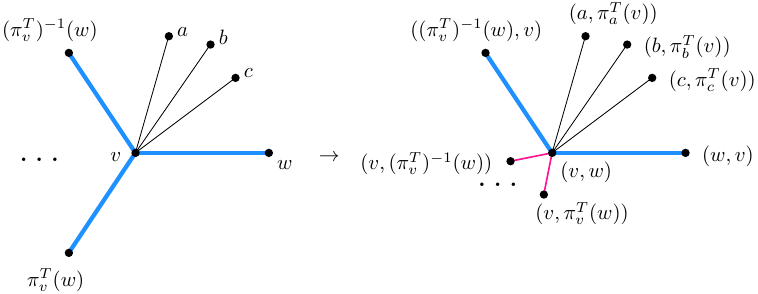}
    \end{center}
    \caption{An illustration of the formal construction of $H''$.}
    \label{fig:orderings}
    \end{figure}

    See \Cref{fig:orderings}.
    It follows from the construction that $H''$ is a planar graph
    and $|V(H'')| = \sum_{v \in V(T)} |N_T(v)| = 2|E(T)| = 2|V(T)|-2 = 2|V(H)|$.
    Let $u'$ be the unique neighbor of $u$ in $T$.
    Let $u_1=u$ and $u_2=u'$.
    For every positive integer $i$, let $u_{i+2} = \pi^T_{u_{i+1}}(u_i)$.
    The walk $u_1,u_2,\ldots,u_{2|E(T)|+1}$ traverses twice through all the edges of $T$, $u_{2|E(T)|+1}=u_1=u$, and $u_{2|E(T)|+2}=u_2=u'$.
    Therefore $(u_1,u_2)(u_2,u_3)\cdots(u_{2|E(T)|},u_1)$ is a Hamiltonian cycle of $H''$.
    Finally, for every $v\in V(H')$, let 
    \begin{align*}
    C_v&=\set{(v,w) \mid w \in N_{T}(v)}. \\
    \intertext{Note that $(C_v\mid v\in V(H'))$ is a model of $H'$ in $H''$.
    Let $F$ be the union of the first row and the first column of $\boxplus_{|V(H'')|}$. Applying~\Cref{lemma:RWT_planar_hamiltonian} to $H''$ and $(u,u')$,  we fix $(B_x \mid x \in V(H''))$ a model of $H''$ in $\boxplus_{|V(H'')|}$ such that $B_{(u,u')}\subseteq F$. Composing the two models together, let}
    D_v &= \textstyle\bigcup_{p\in C_v} B_p\\
    \intertext{for each $v\in V(H')$. Note that $(D_v\mid v\in V(H'))$ is a model of $H'$ in $\boxplus_{|V(H'')|}$. We adapt this model of $H'$ in $\boxplus_{|V(H'')|}$ to an $F$-outer-rooted model of $H$ in $\boxplus_{|V(H'')|}$. 
    Recall that $C_u=\set{(u,u')}$ and therefore $D_u=B_{(u,u')}\subseteq F$. 
    Let $v$ be a vertex of the outer face of~$H$. 
    Since $uv$ is an edge in $H'$, 
    there is an edge between $D_u=B_{(u,u')}$ and $D_v$ in $\boxplus_{|V(H'')|}$.  
    Let this edge be $ab$ with $a \in B_{(u,u')}$. 
    We define}
    E_v &= 
    \begin{cases}
        D_v \cup \{a\}&\textrm{if $D_v\cap F = \emptyset$,}\\
        D_v&\textrm{if $D_v\cap F \neq \emptyset$.}
    \end{cases}
\end{align*}
    Also for each vertex $v$ that is not in the outer face of $H$, we put $E_v=D_v$.
    Finally, we claim that $(E_v\mid v\in V(H))$ is an $F$-outer-rooted model of $H$ in $\boxplus_{|V(H'')|}$. 
    Indeed, note that the sets $(E_v\mid v\in V(H))$ are pairwise disjoint as this is the case for $(D_v\mid v\in V(H))$ and a fixed vertex $a\in F$ can be added to at most one set $E_v$ since every vertex of the outer face of a grid has at most one neighbor which not in the outer face.

    Let $G$ be a graph and $S \subset V(G)$.
    Assume that $G$ contains an $S$-outer-rooted model of $\boxplus_{2|V(H)|}$. 
    Let $F$ be the union of the first row and the first column of $\boxplus_{2|V(H)|}$. 
    We proved above that $\boxplus_{2|V(H)|}$ contains an $F$-outer-rooted model of $H$. 
    Composing the two models we obtain an $S$-outer-rooted model of $H$ in $G$.
    By contraposition, if $G$ has no $S$-outer-rooted model of $H$, then $G$ has no $S$-outer-rooted model of $\boxplus_{2|V(H)|}$, and so, by~\Cref{thm:Stw-grids}, $\tw(G,S) \leq 3 \gm(6144\cdot (2|V(H)|)^4) + 1$.
\end{proof}

\section{Excluding a tangle focused on a set of vertices}
\label{appendix}

Here, we prove~\Cref{thm:Marx_tangles} while closely following the proof of~\cite[statement~1.3]{Marx2017} (which is both more precise and less general than what we need). 
Recall that $\gm$ is the minimum function such that for every positive integer $\ell$, if a graph $G$ has no model of $\boxplus_\ell$, then $\tw(G) \leq \gm(\ell)$.
We need the following result by Kawarabayashi, Thomas, and Wollan~\cite{kawarabayashi2021quicklyexcludingnonplanargraph}.

\begin{theorem}[Lemma 14.6 in~\cite{kawarabayashi2021quicklyexcludingnonplanargraph}]\label{thm:grid_minor_thm_tangles}
    For every positive integer $\ell$, for every graph $G$,
    for every tangle $\mathcal{T}$ of $G$ of order at least $3\gm(6\ell^2)+1$,
    there exists a model $(B_x \mid x \in V(\boxplus_\ell))$ of $\boxplus_\ell$
    in~$G$ such that for every separation $(A,B) \in \mathcal{T}$
    of order at most $\ell-1$, there is no row $R$ of $\boxplus_\ell$ such that
    $\bigcup \{B_x \mid x \in R \} \subseteq V(A)$.
\end{theorem}

Note that this result is stated in~\cite{kawarabayashi2021quicklyexcludingnonplanargraph} in terms of walls instead of grids, but this statement can be deduced using the fact that their $\ell$-walls contain $\boxplus_\ell$ as a minor.

Let $G$ and $H$ be graphs.
A \emph{pseudomodel} of $H$ in $G$ is a family
$(B_x \mid x \in V(H))$ of nonempty subsets of $V(G)$ such that
for every $xy \in E(H)$, there is an edge between $B_x$ and $B_y$ in $G$.

Let $S \subseteq V(G)$, let $G'$ be a grid, and let $k$ be positive integers.
We say that a model $(B_x \mid x \in V(G'))$ of $G'$ in $G$ is \emph{$(S,k)$-augmentable in $G$} if there are $k$ distinct vertices $x_1, \dots, x_k$ of the first row of $G'$
and $k$ pairwise disjoint paths $Q_1, \dots, Q_k$ in $G$
such that for every $i \in [k]$,
$Q_i$ is an $S$--$B_{x_i}$ path and is internally disjoint from $\bigcup \{B_{y} \mid y \in V(G')\}$.

Let $\ell$ be a positive integer.
When $J$ is a subgraph of $\boxplus_\ell$, we define the \emph{boundary} of $J$ in $\boxplus_\ell$ as the set of all vertices $x \in V(J)$ such that there exists $y \in V(\boxplus_\ell)$ adjacent to $x$ in $\boxplus_\ell$ with $xy \notin E(J)$.

\begin{lemma}\label{lemma:induction_rooted_grid}
    Let $k, \ell, n$ be positive integers such that $n \geq (2k+1)(2\ell-1)$.
    Let $G$ be a graph and $S \subseteq V(G)$.
    Let $J$ be a subgraph of $\boxplus_n$ such that $V(J)$ contains at least one row of $\boxplus_n$. 
    Let $\beta$ be a subset of vertices of $J$ containing the boundary of $J$ in $\boxplus_n$.
    Suppose that $|\beta| \leq 2k$.
    Let $(B_x \mid  x \in V(J))$ be a pseudomodel of $J$ in $G$ such that
    \begin{enumerate}
        \item for every $x \in V(J) \setminus \beta$, $B_x$ induces a connected subgraph of $G$, \label{item-i}
        \item for every $x \in \beta$, every component of $G[B_x]$ intersects $S$,\label{item-ii}
        \item there is no separation $(A,B)$ of $G$ of order at most $2k-1$
        and a row $R$ of $\boxplus_n$ contained in $V(J)$ such that
        $S \subseteq V(A)$ and $\bigcup\{B_x \mid x \in R\} \subseteq V(B)$. \label{item-iii}
    \end{enumerate}
    Then, there is a subgraph $H$ of $J$ isomorphic to the $n \times \ell$ grid
    disjoint from $\beta$ and such that $(B_x \mid x \in V(H))$ is 
    $(S,k)$-augmentable in $G$.
\end{lemma}

\begin{proof}
    Consider a counterexample $G,S,J,\beta,(B_x \mid x \in V(J))$
    for which $(|V(G)|+|E(G)|,2k-|\beta|)$ is minimal in the lexicographic order.

    \begin{claim}\label{item-iv}
        There is no separation $(A,B)$ of $G$ of order at most $2k$ with $B \neq G$ and row $R$ of~$\boxplus_n$ contained in $V(J)$ such that $S \subseteq V(A)$ and $\bigcup \{B_x \mid x \in R\} \subseteq V(B)$. 
    \end{claim}

    \begin{proofclaim}
        Suppose for contradiction that there is a separation $(A,B)$ of $G$ of order at most $2k$ with $B \neq G$
        and a row $R$ of $\boxplus_n$ contained in $V(J)$
        such that $S \subseteq V(A)$ and $\bigcup \{B_x \mid x \in R\} \subseteq V(B)$.
        We fix such a separation $(A,B)$ and a row $R$. Note that by \ref{item-iii}, $(A,B)$ has order exactly $2k$.
    
        Let $G' = B$ and $S' = V(A) \cap V(B)$.
        Observe that $|V(G')|+|E(G')| < |V(G)|+|E(G)|$ since $G'=B\neq G$.
        Let $J'$ be the subgraph of $J$ with vertex set $\{x \in V(J) \mid B_x \cap V(B) \neq \emptyset\}$
        and $xy \in E(J)$ is an edge of $J'$ when $x,y\in V(J')$ and either $B_x \cap V(A) = \emptyset$ or $B_y \cap V(A) = \emptyset$.
        Let $\beta' = \{x \in V(J') \mid B_x \cap V(A) \neq \emptyset\}$.
        For every $x \in V(J')$, let $B'_x = B_x \cap V(B)$.
        Then, $(B'_x \mid x \in V(J'))$ is a pseudomodel of $J'$ in $G'=B$.
        We claim that $G',S',J',\beta',(B'_x \mid x \in V(J'))$ satisfy the hypothesis of the lemma.

        First note that $|\beta'|\leq |V(A)\cap V(B)| \leq 2k$.
        Next we argue that the boundary of $J'$ in $\boxplus_n$ is contained in $\beta'$.
        Let $x$ be a vertex of the boundary of $J'$. 
        Suppose that $x\in \beta$.
        By~\ref{item-ii}, $G[B_x]$ intersects $S\subseteq V(A)$ so $x\in \beta'$.
        Now, suppose that $x\notin \beta$. 
        Since $x$ is in the boundary of $J'$ but not in the boundary of $J$,
        $x$ is incident in $J$ to an edge $xy$ not in $J'$.
        Hence, either $y \in V(J')$ and so $x \in \beta'$, or $y\in V(J)-V(J')$.
        In the second case, by definition of $J'$, $B_y \subseteq V(A)\setminus V(B)$ so $B_x\cap V(A)\neq \emptyset$, i.e. $x\in \beta'$.
        This proves that the boundary of $J'$ is contained in $\beta'$.
    
        Since $\bigcup \{B_x \mid x \in R\} \subseteq V(B)$,
        we have $R \subseteq V(J')$, and so, $V(J')$ contains a row of $\boxplus_n$.
        For every $x \in V(J') \setminus \beta'$,
        we have $B_x \cap V(A) = \emptyset$, and since $S \subseteq V(A)$,
        we have $B_x \cap S = \emptyset$.
        Therefore, $B_x = B'_x$ induces a connected subgraph of $G'$, which gives~\ref{item-i}.
        Next, let $x \in \beta'$.
        If $x \in \beta$, then every component of $G[B_x]$ intersects $S\subseteq V(A)$,
        which implies that every component of $G'[B'_x]$ intersects $V(A) \cap V(B) = S'$.
        If $x \notin \beta$, then $G[B_x]$ is connected, and so, by construction, every component of $G'[B'_x]$ intersects $V(A) \cap V(B) = S'$.
        This yields~\ref{item-ii}.
        Finally, suppose to the contrary that there is a separation $(A',B')$ of $G'$ of order at most $2k-1$
        such that $S' \subseteq V(A')$ and $\bigcup \{B_x \mid x \in R'\} \subseteq B'$ for some row $R'$
        of $\boxplus_n$ contained in $J'$.
        Then $(A \cup A',B')$ is a separation of $G$ of order at most $2k-1$
        such that $\bigcup \{B_x \mid x \in R'\} \subseteq B'$ and $S \subseteq V(A \cup A')$,
        which is a contradiction, and hence, we obtain~\ref{item-iii}.
        
        As claimed, $G',S',J',\beta',(B'_x \mid x \in V(J'))$ satisfy the hypothesis of the lemma.
        By minimality of $|V(G)|+|E(G)|$, there exists a subgraph $H$ of $J'$ isomorphic to the $n \times \ell$ grid
        disjoint from $\beta'$ and such that $(B_x' \mid x \in V(H))$ is $(S',k)$-augmentable in $G$.
        In other words, there exist distinct vertices $x_1, \dots, x_k$ of the first row
        of $H$ and $k$ disjoint paths $Q'_1, \dots, Q'_k$ in $G'$
        such that for every $i \in [k]$,
        $Q'_i$ is an $S'$--$B_{x_i}'$ path which is internally disjoint from $\bigcup \{B'_{y} \mid y \in V(H)\}$.
        By Menger's Theorem and \ref{item-iii}, there are $2k = |V(A) \cap V(B)|$ pairwise disjoint $S$--$(V(A) \cap V(B))$ paths in $G$,
        say $P_z$ for $z \in V(A) \cap V(B)$, with $z \in V(P_z)$ for every $z \in V(A) \cap V(B)$.
        For every $i \in [k]$, let $u_i$ be the vertex of $Q'_i$ in $S'$,
        and $Q_i$ the concatenation of $P_{u_i}$ with $Q'_i$.
        Observe that $Q_i$ is a path from $S$ to $B'_{x_i}$ in $G$.
        Since $V(H) \cap \beta' = \emptyset$,
        we have $B'_y = B_y$ for every $y \in V(H)$.
        In particular, $Q_i$ is a path from $S$ to $B_{x_i}$ in $G$ which is
        internally disjoint from $\bigcup\{B_y \mid y \in V(H)\}$.
        Therefore, $(B_x \mid x \in V(H))$ is $(S,k)$-augmentable in $G$, contradicting the fact that $G,S,J,\beta,(B_x \mid x \in V(J))$ is a counterexample.
    \end{proofclaim}

    \begin{claim} \label{claim:singleton}
        For every $x \in V(J) \setminus \beta$, $B_x$ is a singleton.
    \end{claim}
    
    \begin{proofclaim}
        Suppose by contradiction that there exists $x \in V(J) \setminus \beta$ such that $B_x$ is not a singleton.
        Then, by \ref{item-i}, $G[B_x]$ is connected,
        and so there exists an edge $uv \in E(G[B_x])$.
        Let $G'$ be the graph obtained from $G$ by contracting the edge $uv$
        into a single vertex~$w$.
        For every $y \in V(J)$, let
        \begin{align*}
            B'_y &=
            \begin{cases}
                (B_y \setminus \{u,v\}) \cup \{w\} &\textrm{if $y=x$,}\\
                B_y &\textrm{if $y \neq x$,}
            \end{cases} \\
        \intertext{and let}
            S' &=
            \begin{cases}
                (S \setminus\{u,v\}) \cup \{w\} &\textrm{if $u \in S$ or $v \in S$,} \\
                S &\textrm{otherwise.}
            \end{cases}
        \end{align*}
        We claim that $G',S',J,\beta,(B'_y \mid y \in V(J))$ satisfy the hypothesis of the lemma.
        \Cref{item-i} is true since contracting an edge in a branch set does not change its connectivity.
        For every $x\in \beta$, since every component of $G[B_x]$ intersects $S$, every component of $G'[B'_x]$ intersects $S'$.
        Therefore, \cref{item-ii} holds.
        We now prove \cref{item-iii}.
        Suppose to the contrary that there is a separation $(A',B')$ of $G'$ of order at most $2k-1$
        and a row $R$ of $\boxplus_n$ contained in $V(J)$ such that $S' \subseteq V(A')$ and $\bigcup \{B'_y \mid y\in R\} \subseteq V(B')$.
        We define the separation $(A,B)$ of $G$ by
        \begin{align*}
            A &= 
            \begin{cases}
                G[(V(A') \setminus \{w\}) \cup \{u,v\}] &\textrm{if $w \in V(A')$,} \\
                G[V(A')] &\textrm{otherwise,}
            \end{cases} \\
            B &= 
            \begin{cases}
                G[(V(B') \setminus \{w\}) \cup \{u,v\}] \setminus E(A) &\textrm{if $w \in V(B')$,} \\
                G[V(B')] \setminus E(A) &\textrm{otherwise.}
            \end{cases}
        \end{align*}
        Observe that $(A,B)$ is a separation of $G$ of order at most $2k-1-1+2 = 2k$ such that $S \subseteq V(A)$ and $\bigcup \{B_y \mid y \in R\} \subseteq V(B)$.
        If $G \neq B$, this contradicts \Cref{item-iv}.

        Now, suppose $G=B$.
        If $|S| = |S'|$, then
        $|S| \leq |V(A')| = |V(A')\cap V(B')| \leq 2k-1$.
        It follows that $(G[S]-E(G[S]),G)$ is a separation of $G$ of order at most $2k-1$ which contradicts \ref{item-iii} for $G,S,J,\beta,(B_y \mid y \in V(J))$.
        If $|S| = |S'|+1$, then $\{u,v\}\subseteq S$ and $|S|\leq 2k$.
        Since $\{u,v\}\subseteq B_x$ and $x\notin \beta$, by~\ref{item-ii}, $|\beta| \leq |S|-2 \leq 2k-2$.
        Then, replacing $\beta$ by $\beta \cup \{x\}$ yields another counterexample, which contradicts the minimality of $2k-|\beta|$.
        
        Therefore,~\ref{item-iii} holds for $G',S',J,\beta,(B'_y \mid y \in V(J))$, and so, by minimality of $|V(G)|+|E(G)|$,
        there is a subgraph $H$ of $J$ isomorphic to the $n \times \ell$ grid disjoint from $\beta$ and such that $(B_y' \mid y \in V(H))$ is $(S',k)$-augmentable in $G$.
        By construction, $H$ is disjoint from $\beta$ and $(B_y \mid y \in V(H))$ is $(S,k)$-augmentable in $G$ by uncontracting $uv$, which contradicts the fact that $G,S,J,\beta,(B_x \mid x \in V(J))$ is a counterexample.
    \end{proofclaim}

    \begin{claim}\label{claim:no_path_outside}
        For every $x, y \in V(J)$ with $x\neq y$ and $xy\not\in E(J)$, there is no $B_x$--$B_y$ path in $G$ internally disjoint from $\bigcup\{B_z \mid z \in V(J)\}$.
    \end{claim}
    \begin{proofclaim}
        Suppose by contradiction that there exists $x,y \in V(J)$ distinct and non adjacent in $J$
        and an $B_x$--$B_y$ path $P$
        in $G$ internally disjoint from $\bigcup\{B_z \mid z \in V(J)\}$.
        Then, the first edge $uv$ along $P$ is such that $u \in B_x$ and $v \not\in \bigcup\{B_x \mid z \in N_J[x]\}$.
        Let $G' = G \setminus \{uv\}$.
        We claim that $G',S,J,\beta,(B_y \mid y \in V(J))$ satisfies the hypothesis of the lemma.
        Since $uv$ is not contained in any branch set, removing it does not affect~\ref{item-i} and~\ref{item-ii}.
        For~\ref{item-iii}, suppose to the contrary that there is a separation $(A',B')$ of $G'$ of order at most $2k-1$
        and a row $R$ of $\boxplus_n$ contained in $V(J)$ such that $S \subseteq V(A)$ and $\bigcup \{B_y \mid y \in R\} \subseteq V(B)$.
        By item \ref{item-iii} for $G,S,J,\beta,(B_y \mid y \in V(J))$, the pair $(G[V(A')],G[V(B')]\setminus E(G[V(A')]))$ is not a separation of $G$.
        Therefore, without loss of generality, assume that $u \in V(A') \setminus V(B')$
        and $v \in V(B') \setminus V(A')$.
        Then consider the separation $(A,B)$ of $G$ defined by
        \begin{align*}
            A &= G[V(A') \cup \{v\}], \\
            B &= G[V(B')] \setminus E(A).
        \end{align*}
        Observe that $(A,B)$ is a separation of $G$ of order at most $2k-1+1 = 2k$
        such that $S \subseteq V(A)$ and $\bigcup \{B_y \mid y \in R\} \subseteq V(B)$.
        If $G \neq B$, this contradicts \Cref{item-iv}.
        If $G = B$, then $V(A) = V(A) \cap V(B)$ has size at most $2k$.
        Since $v \not\in V(A')$, we have $v \not \in S$ and so $S \subseteq V(A) \setminus \{v\}$,
        and it follows that $|S| \leq 2k-1$.
        Then, the separation $(G[S] \setminus E(G[S]), G)$ contradicts \ref{item-iii} for $G,S,J,\beta,(B_y \mid y \in V(J))$.
    \end{proofclaim}
    
    Since $|\beta| \leq 2k$ and $n \geq (2k+1)(2\ell-1)$, by the pigeonhole principle,
    there are $2\ell-1$ consecutive rows $R_1, \dots, R_{2\ell-1}$ of $\boxplus_n$
    which are disjoint from $\beta$.
    This implies that either $\bigcup \{R_i \mid i \in [2\ell-1] \} \subseteq V(J)-\beta$,
    or $\bigcup \{ R_i \mid i \in [2\ell-1] \} \cap V(J) = \emptyset$.
    
    First, suppose that $\bigcup \{R_i \mid i \in [2\ell-1] \} \cap V(J) = \emptyset$.
    By assumption, there is a row $R$ of $\boxplus_n$ contained in $V(J)$.
    Since every column of $\boxplus_n$ contains a vertex in $\beta$, $|\beta| \geq n > 2k$, a contradiction.
    
    Now, suppose that $\bigcup \{ R_i \mid i \in [2\ell-1]\} \subseteq V(J)\setminus\beta$.
    By Menger's Theorem, either there is a separation $(A,B)$ of $G$ of order at most $2k-1$ with $S \subseteq V(A)$ and $\bigcup \{B_x \mid x\in R_\ell\} \subseteq V(B)$,
    or there are $2k$ pairwise disjoint $S$--$\bigcup \{B_x \mid x\in R_\ell\}$ paths $Q_1,Q_2,\dots, Q_{2k}$ in $G$.
    By~\ref{item-iii}, the latter holds.
    Let $j \in [2k]$.
    We order the vertices in $Q_j$ from the endpoint of $Q_j$ in $S$ to the endpoint in $\bigcup\{B_x \mid x \in R_\ell\}$.
    Let $u_j$ be the first vertex of $Q_j$ in $\bigcup \{B_x \mid x\in \bigcup_{i \in [2\ell-1]} R_i\}$ and 
    let $y_j \in R_\ell$ be such that the first vertex along $Q_j$ in $\bigcup\{B_x \mid x \in R_\ell\}$ belongs to $B_{y_j}$.
    By \Cref{claim:singleton}, $B_{y_1},B_{y_2},\dots,B_{y_{2k}}$ are singletons because $\{y_j \mid j\in [2k]\}\subseteq V(J)-\beta$.
    Therefore, $y_1,y_2,\dots,y_{2k}$ are distinct vertices.
    There exists $X \subseteq [2k]$ of size at least $k$
    such that either $\{u_j \mid j\in X\} \subseteq \{ B_x \mid x \in \bigcup_{i=1}^{\ell} R_i\}$
    or $\{u_j \mid j\in X\} \subseteq \bigcup \{ B_x \mid x \in \bigcup_{i=\ell}^{2\ell-1} R_i\}$.
    Without loss of generality, assume the former holds.
    Let $H$ be the subgraph of $\boxplus_n$ induced by $\bigcup_{i=1}^{\ell} R_i$.
    Observe that $H$ is isomorphic to the $n \times \ell$ grid,
    and the vertices $y_j$ for each $j \in X$ are in the first row of $H$, i.e the row corresponding to $R_\ell$.
    By \Cref{claim:no_path_outside}, the paths $Q_j$ for each $j \in X$ are internally disjoint from $V(H)$ since $\bigcup \{B_x\mid x \in R_\ell\}$ intersects every path between $V(H)$ and $\bigcup \{B_x\mid x\in \bigcup_{i=\ell}^{2\ell-1} R_i \}$.
    It follows that $(B_x \mid x \in V(H))$ is $(S,k)$-augmentable in $G$, which ends the proof.
\end{proof}

We can now prove \Cref{thm:Marx_tangles}.

\begin{figure}[!htbp]
    \centering 
    \includegraphics{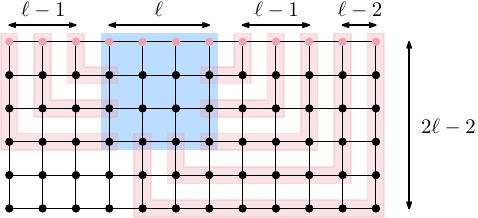} 
    \caption{
        From a model of a $(4\ell-4)\times (2\ell-2)$ grid with each branch set in the first row intersecting $S$, it is easy to construct an $S$-outer-rooted model of $\boxplus_\ell$.
    } \label{fig:grids}
\end{figure}

\begin{proof}[Proof of \Cref{thm:Marx_tangles}]
    Let $\ell$ be a positive integer and set $k = 4\ell-4$, $\ell' = 2\ell-2$ and $n = (2k+1)(2\ell'-1)$.
    Recall that by~\cite{Chuzhoy2021}, $\gm$ is upper bounded by a polynomial function, and $\gm$ is non-decreasing.
    Let $K = 3\gm(6144\ell^4) + 1 \geq 3\gm(6n^2)+1$.
    Let $G$ be a graph and $S \subset V(G)$.
    We assume that $G$ has a tangle $\mathcal{T}$ of $(G,S)$ of order at least $K$.
    By \Cref{thm:grid_minor_thm_tangles},
    there is a model $(B_y \mid y \in V(\boxplus_n))$ of $\boxplus_n$
    in $G$ such that for every separation $(A,B) \in \mathcal{T}$,
    of order at most $n-1$, there is no row $R$ of $\boxplus_n$ such that $\bigcup\{ B_y \mid y \in R\} \subseteq V(A)$.
    Now, for every separation $(A,B)$ of $G$ of order at most $2k-1$,
    either $(A,B) \in \mathcal{T}$ and then $S \not\subseteq V(A)$,
    or $(B,A) \in \mathcal{T}$ and then there is no row $R$ of $\boxplus_n$ such that $\bigcup\{ B_y \mid y \in R\} \subseteq V(B)$.
    Therefore, by \Cref{lemma:induction_rooted_grid} applied
    for $J = \boxplus_n$ and $\beta = \emptyset$,
    there is $H \subseteq \boxplus_n$ isomorphic to the $\ell' \times n$ grid
    such that $(B_y \mid y \in V(H))$ is $(S,k)$-augmentable in $G$.
    In other words, there exists $x_1, \dots, x_k$ distinct vertices in the first row of $H$,
    and $k$ disjoint paths $Q_1, \dots, Q_k$ where for every $i\in [k]$, $Q_i$ is an $S$--$B_{x_i}$ path
    internally disjoint from $\bigcup\{ B_y \mid y \in V(H)\}$.
    By respectively adding vertices of $Q_1,\dots,Q_k$ to the branch sets of $x_1,\dots,x_k$,
    we obtain a model of the $\ell' \times n$ grid in $G$ such that each of the branch sets corresponding to $x_1, \dots, x_k$ intersects $S$.
    We denote by $(C_y \mid y \in [\ell'] \times [n])$ this model.
    Next, we will ``contract some horizontal edges'' on columns not containing an $x_i$ for some $i\in [k]$.
    Formally, let $n_0 = \min \{j\in [n] \mid \{x_1,\dots,x_k\}\subseteq\{(1,1),\dots,(1,j)\}\}$. 
    For every $j\in [n_0]$, let $x(j) = \min\{j'\geq j \mid (1,j')\in \{x_1,\dots,x_k\}\}$.
    For every $j\in [k]$ and every $i\in [\ell']$, let $D_{(i,j)} = \bigcup \{ C_{(i,j')} \mid j'\in [n_0], x(j')=j \}$.
    As a result, $(D_y \mid y\in [\ell'] \times [k])$ is a model of the $\ell' \times k = (2\ell-2) \times (4\ell-4)$ grid in $G$ 
    such that the branch set of every vertex of the first row intersects $S$.
    Finally, one can easily find $4\ell-4$ disjoint paths from the first row of a $(2\ell-2) \times (4\ell-4)$ grid to the boundary of a subgraph $H'$ isomorphic to $\boxplus_\ell$ as shown in~\Cref{fig:grids}.
    Adding vertices of these paths to the branch sets of their endpoints in $H'$ gives an $S$-outer-rooted model of $\boxplus_\ell$ in $G$.
    This completes the proof.
\end{proof}

\end{document}